\documentclass[12pt]{article}
\usepackage{amsmath,amssymb,amsfonts}
\usepackage{epsfig}
\usepackage{tikz}
\usepackage{graphics}
\usepackage{bbm}
\usepackage{bm, bbm}
\usepackage{todonotes}
\usepackage[thmmarks]{ntheorem}
\usepackage{ifthen}
\usepackage{listings}
\usepackage{bbm}
\usepackage{hyperref}
 \hypersetup{colorlinks=true, citecolor=black,%
                 filecolor=black,%
                 linkcolor=black,%
                 urlcolor=black,  pdftex}

\usetikzlibrary[graphs, arrows, backgrounds, intersections, positioning, fit, petri, calc, shapes, decorations.pathmorphing]
\usepgflibrary{patterns}

\definecolor{darkred}{RGB}{105,0,0}

\makeatletter
\newtheoremstyle{prime}%
 {\item[\hskip\labelsep \theorem@headerfont ##1\ \theorem@separator]}%
{\item[\hskip\labelsep \theorem@headerfont ##1\ ##3' \theorem@separator]}
\makeatother

\makeatletter
\newtheoremstyle{proofof}
{\item[\hskip\labelsep \theorem@headerfont ##1\ \theorem@separator]}%
{\item[\hskip\labelsep \theorem@headerfont ##1\ ##3\theorem@separator]}
\makeatother

\newtheorem{theorem}{Theorem}

\newtheorem{proposition}[theorem]{Proposition}

\newtheorem{corollary}[theorem]{Corollary}

\newtheorem{claim}{Claim}

\theoremstyle{prime}
\newtheorem{stheorem}{Theorem}

\def \QD1 {\hfill $\spadesuit$}

\newcommand{\case}[2]{\noindent {\bf Case #1\/:} {\it #2}}

\newcommand{\DF}[1]{{\bf #1\/}}

\newcommand{\set}[2]{\{#1 \;|\; #2 \}}
\newcommand{\ems}{\varnothing}
\newcommand{\sm}{\setminus}

\newcommand{\De}{\Delta}
\newcommand{\de}{\delta}

\newcommand{\cn}{\chi}
\newcommand{\lcn}{\chi^{\ell}}

\newcommand{\f}{\varphi}

\newcommand{\NP}  {{\sf NP}}
\newcommand{\NPC} {{\NP}-complete}

\newcommand{\cB}{{\cal B}}

\newcommand{\cV}{{\cal V}}

\newcommand{\col}{{\rm col}}

\newcommand{\nato}{\mathbb{N}_0}

\newcommand{\sref}[1]{\ref{#1}'}

\numberwithin{equation}{section}

\theoremstyle {nonumberplain}
\theoremseparator{:}
\theorembodyfont{\normalfont}
\theoremsymbol{\rule{1ex}{1ex}}
\newtheorem{proof}{Proof}

\theoremstyle{proofof}
\newtheorem{proofof}{Proof of}

\theoremsymbol{$\square$}
\newtheorem{proof2}{Proof}

\begin{document}
\title{\bf Partitions of hypergraphs under variable degeneracy constraints}

\author{{{
Thomas Schweser}\thanks{partially supported by DAAD, Germany (as part of BMBF) and by the Ministry of Education Science, Research and Sport of the Slovak Republic within the project 57320575}
\thanks{
Technische Universit\"at Ilmenau, Inst. of Math., PF 100565, D-98684 Ilmenau, Germany. E-mail
address: thomas.schweser@tu-ilmenau.de}}
\and{{Michael Stiebitz}\footnotemark[1]
\thanks{
Technische Universit\"at Ilmenau, Inst. of Math., PF 100565, D-98684 Ilmenau, Germany. E-mail
address: michael.stiebitz@tu-ilmenau.de}}}

\date{}
\maketitle

\begin{abstract}
\noindent The paper deals with partitions of hypergraphs into induced subhypergraphs satisfying constraints on their degeneracy. Our hypergraphs may have multiple edges, but no loops. Given a hypergraph $H$ and a sequence $f=(f_1,f_2, \ldots, f_p)$ of $p\geq 1$ vertex functions
$f_i:V(H) \to \nato$ such that $f_1(v)+f_2(v)+ \cdots + f_p(v)\geq d_H(v)$ for all $v\in V(H)$, we want to find a sequence $(H_1,H_2, \ldots, H_p)$ of vertex disjoint induced subhypergraphs containing all vertices of $H$ such that each hypergraph $H_i$ is strictly $f_i$-degenerate, that is, for every non-empty subhypergraph $H'\subseteq H_i$ there is a vertex $v\in V(H')$ such that $d_{H'}(v)<f_i(v)$. Our main result in this paper says that such a sequence of hypergraphs exists if and only if $(H,f)$ is not a so-called hard pair. Hard pairs form a recursively defined family of configurations, obtained from three basic types of configurations by the operation of merging a vertex. Our main result has several interesting applications related to generalized hypergraph coloring problems.
\end{abstract}



\section{Introduction and main results}

Motivated by a question due to Paul Erd\H{o}s, Lov\'asz \cite{Lovasz66} proved the first result concerning partitions of simple graphs into parts with bounded maximum degree.
A partition of a hypergraph is a sequence of induced subhypergraphs (possibly empty) such that each vertex belongs to exactly one part of the partition. In this paper, we examine partitions of hypergraphs with respect to the degeneracy of the hypergraph. Therefore, we say that, given a hypergraph $H$ and a function $f:V(H) \to \mathbb{N}_{0}$, a hypergraph $H$ is strictly $f$-degenerate if for any non-empty subhypergraph $G$ of $H$ there exists a vertex $v \in V(G)$ such that $d_G(v) < f(v)$. In this paper, our aim is to find sufficient and necessary conditions for a hypergraph $H$ to admit a partition such that each part $H_i$ of the partition is strictly $f_i$-degenerate for some function $f_i$. An analogue problem for the class of simple graphs was considered by Borodin, Kostochka and Toft \cite{BorKosToft} in 2000. We show how to extend their interesting theorem to the class of hypergraphs and also give several applications.

\subsection{Hypergraphs, basic concepts}

A \DF{hypergraph} $H=(V,E,i)$ is a triple consisting of two finite sets,
$V$ and $E$, and a function $i$ from $E$ to the power set $2^V$, such that $|i(e)|\geq 2$ for $e\in E$.
The set $V=V(H)$
is the \DF{vertex set} of
$H$ and its elements are the \DF{vertices} of $H$. The set $E=E(H)$ is the \DF{edge set} of $H$ and its elements are the \DF{edges} of $H$. The function $i=i_H$ is the \DF{incidence function} of $H$ and $i(e)$ is the set of vertices that are \DF{incident} to the edge $e$ in $H$.  A hypergraph $H$ is \DF{empty} if $V(H)=E(H)=\ems$; in this case we write $H=\ems$.

Let $H=(V,E,i)$ be a hypergraph. The number of vertices of $H$ is its \DF{order}, written $|H|$. An edge $e$ with $|i(e)|\geq 3$ is called a \DF{hyperedge}, and an edge $e$ with $|i(e)|=2$ is called an \DF{ordinary edge}. Two distinct edges $e,e'\in E$ with $i(e)=i(e')$ are called \DF{parallel edges}. We say that $H$ is a \DF{simple hypergraph} if $H$ has no parallel edges.  If $|i(e)|=q$ for all $e\in E$, then $H$ is said to be \DF{$q$-uniform}. Thus, a \DF{graph} is a 2-uniform hypergraph, that is, a hypergraph in which each edge is ordinary. Note that in our terminology a graph may have parallel edges, otherwise the graph is said to be a \DF{simple graph}.

A hypergraph $H'$ is a \DF{subhypergraph} of $H$, written $H'\subseteq H$, if $V(H')\subseteq V(H)$, $E(H')\subseteq E(H)$, and $i_{H'}=i_H|_{E(H')}$. If $H'\subseteq H$ and $H'\not =H$, then $H'$ is said to be a \DF{proper subhypergraph} of $H$. For two subhypergraphs $H_1$ and $H_2$ of $H$ we define the \textbf{union} and the \textbf{intersection} of $H_1$ and $H_2$ in the usual way. So $H'=H_1\cup H_2$ is the subhypergraph of $H$ with $V(H')=V(H_1) \cup V(H_2)$, $E(H')=E(H_1)\cup E(H_2)$, and $i_{H'}=i_H|_{E(H')}$; and $H'=H_1\cap H_2$ is the subhypergraph of $H$ with $V(H')=V(H_1) \cap V(H_2)$, $E(H')=E(H_1)\cap E(H_2)$, and $i_{H'}=i_H|_{E(H')}$.

Let $H^1$ and $H^2$ be two \DF{disjoint hypergraphs}, that is, $V(H^1)\cap V(H^2)=\ems$ and $E(H^1) \cap E(H^2)=\ems$. Furthermore, for $i\in {1,2}$, let $v^i\in V(H^i)$ and let
$v^*\not\in (V(H^1) \cup V(H^2))\sm \{v^1,v^2\}$. We obtain a new hypergraph $H$ with $V(H)=((V(H^1) \cup V(H^1))\sm \{v^1,v^2\}) \cup \{v^*\}$, $E(H)=E(H^1) \cup E(H^2)$, and
$$i_H(e)=
\begin{cases}
i_{H^j}(e) & \text{if } e\in E(H^j), v^j \not \in i_{H^j}(e)~(j \in \{1,2\}),\\
(i_{H^j}(e)\sm \{v^j\}) \cup \{v^*\} & \text{if } e\in  E(H^j), v^j \in i_{H^j}(e)~(j \in \{1,2\}).
\end{cases}$$
In this case we say that $H$ is obtained from $H^1$ and $H^2$ by \DF{merging} $v^1$ and $v^2$ to $v^*$.

With a hypergraph $H$ and a vertex set $X\subseteq V(H)$ we associate two new hypergraphs, both having $X$ as its vertex set. First, by $H[X]$ we denote the subhypergraph of $H$ satisfying
$$V(H[X])=X, E(H[X])=\set{e\in E}{i_H(e)\subseteq X} \mbox{, and } i_{H[X]}=i_H|_{E(H[X])}.$$
We call $H[X]$ the \DF{subhypergraph of $H$ induced} by $X$. A hypergraph $H'$ is called an \DF{induced subhypergraph} of $H$ if $V(H') \subseteq V(H)$ and $H'=H[V(H')]$. Secondly, by $H(X)$ we denote the hypergraph satisfying
$$V(H(X))=X, E(H(X))=\set{e\in E}{|i(e)\cap X|\geq 2},$$
and
$$i_{H(X)}(e)=i_H(e) \cap X \mbox{ for all } e\in E(H(X)). $$
We call $H(X)$ the hypergraph obtained by \DF{shrinking $H$ to $X$}. Furthermore, we define $H-X=H[V(H)\sm X]$ and
$H \div X=H(V(H)\sm X)$. When $X=\{v\}$ is a singleton, we denote $H-X$ by $H-v$ and $H\div X$ by $H \div v$. Moreover, if $H'$ is a proper induced subhypergraph of $H$ and if $v \in V(H) \setminus V(H')$, let $H' + v = H[V(H') \cup \{v\}]$. Note that it clearly holds $(H \div u) \div v = (H \div v) \div u$ for all vertices $u \neq v$ from $V(H)$.

Let $H$ be a hypergraph. A vertex set $X\subseteq V(H)$ is called an \DF{independent set} of $H$ if the hypergraph $H[X]$ has no edge; it is called a \DF{clique} of $H$ if ${X \choose 2} \subseteq \set{i(e)}{e\in E(H[X])}$.
We call $H$ a \DF{complete $q$-uniform hypergraph}, where $q\geq 2$ is an integer, if $H$ is a simple hypergraph such that $\set{i(e)}{e\in H}={V(H) \choose q}$. If $H$ is a complete $q$-uniform hypergraph of order $n$, we write $H=K_n^q$. Note that the hypergraph $K_n^n$ with $n\geq 2$ has exactly one edge. For the complete graph $K_n^2$ we also write $K_n$. We write $H=C_n$ for an ordinary cycle as a 2-uniform simple hypergraph of order $n$. A cycle is called \DF{odd} or \DF{even} depending on whether its order is odd or even. For a simple hypergraph $H$ and an integer $t\geq 1$, we denote by $H'=tH$ the hypergraph obtained from $H$ by replacing each edge of $H$ by $t$ parallel edges.

A non-empty hypergraph $H$ is called \DF{connected} if for every vertex set $X$ with $\ems \not= X \varsubsetneq V(H)$ there is at least one edge $e\in E(H)$ such that $i_H(e)$ contains a vertex of $X$ as well as a vertex of $V(G)\sm X$. Equivalently, $H$ is connected if and only if there is a hyperpath in $H$ between any two of its vertices. A \DF{hyperpath of length} $q$ in $H$ is a sequence
$(v_1, e_1, v_2, e_2, \ldots, v_q, e_q, v_{q+1})$ of distinct vertices
$v_1, v_2, \ldots, v_{q+1}$ of $H$ and distinct edges
$e_1, e_2, \ldots, e_q$ of $H$ such that
$\{v_i,v_{i+1}\}\subseteq i_H(e_i)$ for $i=1, 2, \ldots, q$. A (connected) \DF{component} of a nonempty hypergraph $H$ is a maximal connected subhypergraph.

A vertex $v$ of a hypergraph $H$ is called a \DF{separating vertex} of $H$ if $H$ is the union of two induced subhypergraphs $H_1$ and $H_2$ satisfying $V(H_1)\cap V(H_2)=\{v\}$ and $|H_i|\geq 2$ for $i\in \{1,2\}$. It is easy to show that if $H$ is a connected hypergraph, then a vertex $v$ of $H$ is non-separating if and only if the hypergraph $H\div v$ is empty or connected.

A \DF{block} of a hypergraph $H$ is a maximal connected subhypergraph of $H$
that has no separating vertex. Let $\cB(H)$ denote the set of all blocks of $H$ and, given a vertex $v \in V(H)$, let $\mathcal{B}_v(H)$ be the set of all blocks from $\cB(H)$ that contain $v$. Note that $\cB(\ems)=\ems$ and that every block of a non-empty hypergraph $H$ is a connected induced subhypergraph of $H$. As for graphs it is not difficult to show that any two distinct blocks of a hypergraph $H$ have at most one vertex in common, and a vertex of $H$ is a separating vertex of $H$ if and only if it belongs to more than one block. A block of $H$, which contains at most one separating vertex of $H$, is called an \DF{end-block} of $H$. If $H$ contains a separating vertex, then $H$ has at least two end-blocks.

\subsection{Degeneracy of hypergraphs}

Let $H$ be a hypergraph. A vertex $v\in V(H)$ is \DF{incident} with an edge $e\in E(H)$ if $v\in i_H(e)$. Moreover, two vertices $u \neq v \in V(H)$ are \textbf{adjacent}, if there is an edge $e \in E(H)$ such that $\{u,v\} \subseteq i_H(e)$. In this case we say that $u$ is a \textbf{neighbor} of $v$ and vice versa. For a vertex $v$ of $H$, let
$$E_H(v)=\set{e\in E(H)}{v\in i_H(e)}.$$
The \DF{degree} of $v$ in $H$ is $d_H(v)=|E_H(v)|$. The hypergraph $H$ is said to be \DF{regular} and \DF{$r$-regular} if each vertex of $H$ has degree $r$ in $H$.
As usual, $\de(H)=\min_{v\in V(H)} d_H(v)$ is the \DF{minimum degree} of $H$ and $\De(H)=\max_{v\in V(H)}d_H(v)$ is the \DF{maximum degree} of $H$. If $H=\ems$, then we define $\de(H)=\De(H)=0$. For an ordinary edge $e$ of $H$ with $i_H(e)=\{u,v\}$, we also write $e=uv$ and $e=vu$. For two distinct vertices $u$ and $v$ of $H$, let
$$\mu_H(u,v)=|\set{e\in E(H)}{e=uv}|$$
be the \DF{multiplicity} of $(u,v)$ in $H$. Note that if $v\in V(H)$, then every vertex $u\in V(H)\sm \{u\}$ satisfies
\begin{align}
\label{Equ:degree(uinH-v)}
d_{H\div v}(u)=d_H(u)-\mu_H(u,v).
\end{align}

The hypergraph $H$ is called \DF{strictly $k$-degenerate}, where $k$ is a non-negative integer, if every non-empty subhypergraph $H'$ of $H$ contains a vertex $v$ with $d_{H'}(v)<k$. So $H$ is strictly $0$-degenerate if and only if $H=\ems$, and $H$ is strictly $1$-degenerate if and only if $E(H)=\ems$. The strictly $2$-degenerate graphs are precisely the forests. Lick and White \cite{LickWhite70} defined a graph to be \DF{$k$-degenerate} if each of its non-empty subgraphs has a vertex of degree at most $k$. Thus, being strictly $k$-degenerate is equivalent to being $(k-1)$-degenerate. We have made this deviation from their terminology to express our results in a more natural way. The smallest $k$ for which the hypergraph $H$ is strictly $k$-degenerate is called the \DF{coloring number} of $H$, denoted by $\col(H)$. That the coloring number of simple graphs can be computed by a polynomial time algorithm was observed independently by various researchers including Finck and Sachs \cite{FinckS69}, Matula \cite{Matul68}, and possibly others. Their arguments can easily be extended to hypergraphs.

Let $h$ be a function from $V(H)$ to the set $\nato$ of non-negative integers. The hypergraph $H$ is said to be \DF{strictly $h$-degenerate} if every non-empty subhypergraph $H'$ of $H$ contains a vertex $v$ such that $d_{H'}(v)<h(v)$.

\subsection{Partitions and colorings of hypergraphs}

Let $H$ be an arbitrary hypergraph. A sequence $(H_1,H_2, \ldots, H_p)$ of $p\geq 1$ pairwise vertex disjoint induced subhypergraphs of $H$ such that $V(H)=V(H_1) \cup V(H_2) \cup \cdots \cup V(H_p)$ is called a \DF{partition} and a \DF{$p$-partition} of $H$. Note that in a partition also empty subhypergraphs are allowed.

A \DF{coloring} of $H$ with color set $C$ is a function $\f:V(H) \to C$. Let $\f$ be a coloring of $H$ with color set $C$. If $|C|=k$, we also say that $\f$ is a $k$-coloring of $H$. Furthermore, for each \DF{color} $c\in C$, the set
$$\f^{-1}(c)=\set{v\in V(H)}{\f(v)=c}$$
is called a \DF{color class} of $H$ with respect to $\f$. In many coloring problems it is required to choose the color for each vertex $v$ from an individual list $L(v)$ of available colors. Therefore, we call a function $L:V(H)\to 2^C$ a \DF{list-assignment} of $H$ with color set $C$. For a given list-assignment $L$ of $H$ with color set $C$ and a coloring $\f$ of $H$ with color set $C$, we say that $\f$ is an $L$-coloring of $H$ if $\f(v)\in L(v)$ for all $v\in V(H)$.

The concepts of hypergraph partitions and hypergraph coloring are closely related to each other. Let $H$ be a hypergraph and let $C=\{1,2, \ldots, p\}$ be a color set. If $\f$ is a coloring of $H$ with color set $C$, then $(H_1,H_2, \ldots, H_p)$ with $H_c=H[\f^{-1}(c)]$ for $c\in C$ is a partition of $H$. Conversely, if $(H_1,H_2, \ldots, H_p)$ is a partition of $H$, then the function $\f$ with $\f(v)=c$ if $v\in V(H_c)$ is a coloring of $H$ with color set $C$.

Colorings and partitions of hypergraphs become a subject of interest only when some restrictions
to the color classes, respectively to the parts of the partition, are imposed. For instance,
a coloring or list-coloring of a hypergraph $H$ with color set $C$ is called a \DF{proper coloring},
respectively a \DF{proper list-coloring} of $H$ if each color class is an independent set of $H$ and
induces therefore an edgeless subhypergraph of $H$. Note that in a proper coloring, each color class induces a strictly $1$-degenerate subhypergraph. The \DF{chromatic number} of a hypergraph $H$, denoted by $\cn(H)$, is the least integer $k$ such that
$H$ admits a proper $k$-coloring. Similar, the \DF{list-chromatic number} of $H$, denoted by
$\lcn(H)$, is the least integer $k$ such that $H$ admits a proper $L$-coloring for each list-assignment $L$ satisfying $|L(v)| \geq k$ for all $v \in V(H)$. Since $\lcn(H) = k$ implies that $H$ has a proper $L$-coloring for the list-assignment $L$ with $L(v) = \{1,2, \ldots, k\}$, we obtain
$\cn(H)\leq \lcn(H)$. Furthermore, a simple sequential coloring argument shows that
$$\cn(H)\leq \lcn(H) \leq \col(H) \leq \De(H)+1.$$
Note that the chromatic number and the list chromatic number of a hypergraph $H$ is equal to the chromatic number, respectively list chromatic number of its underlying simple hypergraph, that is, the hypergraph obtained from $H$ be replacing all parallel edges by a single edge.

\subsection{Problem and main result}

Let $H$ be an arbitrary hypergraph. A function $f:V(H) \to \nato^p$ is called a \DF{vector function} of $H$. By $f_i$ we name the $i$th coordinate of $f$, i.e., $f=(f_1,f_2, \ldots, f_p)$. The set of all vector functions of $H$ with $p$ coordinates is denoted by $\cV_p(H)$. For $f\in \cV_p(H)$, an \DF{$f$-partition} of $H$ is a $p$-partiton $(H_1,H_2, \ldots, H_p)$ of $H$ such that $H_i$ is strictly $f_i$-degenerate for all $i\in \{1,2, \ldots, p\}$. If the hypergraph $H$ admits an $f$-partition, then $H$ is said to be \DF{$f$-partitionable}. Note that for $p=1$, the following proposition clearly holds. This observation will be used frequently in the next section.

\begin{proposition} \label{prop_reg-f-part}
Let $H$ be a connected hypergraph, and let $h \in \mathcal{V}_1(H)$. If $d_H(v)=h(v)$ for all $v \in V(H)$, then each proper subhypergraph of $H$ is strictly $h$-degenerate.
\end{proposition}

In the following we will examine, which conditions are sufficient for a hypergraph $H$ in order to admit an $f$-partition. A first natural suggestion would be, given a hypergraph $H$ and a function $f \in \mathcal{V}_p(H)$, that the requirement $f_1(v) + f_2(v) + \ldots + f_p(v) \geq d_H(v)$ for all $v \in V(H)$ is adequate. However, it is not hard to find a huge number of pairs $(H,f)$ for which this condition is not sufficient. The good news is that all of those 'bad pairs' can be characterized nicely. To this end, we introduce the following, recursively defined class of configurations.

Let $H$ be a connected hypergraph and let $f \in \mathcal{V}_p(H)$ be a vector-function for some $p \geq 1$. We say that $H$ is $f$\textbf{-hard}, or, equivalently, that $(H,f)$ is a \textbf{hard pair}, if one of the following four conditions hold.

\begin{itemize}
\item[(1)] $H$ is a block and there exists an index $j \in \{1,2,\ldots,p\}$ such that
$$f_i(v)=
\begin{cases}
d_H(v) & \text{if i=j,}\\
0 & \text{otherwise}
\end{cases}$$
for all $i \in \{1,2,\ldots,p\}$ and for each $v \in V(H)$. In this case, we say that $H$ is a \textbf{monoblock} or a block of type \textbf{(M)}.
\item[(2)] $H=tK_n$ for some $t \geq 1, n \geq 3$ and there are integers $n_1,n_2,\ldots,n_p \geq 0$ with at least two $n_i$ different from zero such that $n_1 + n_2 + \ldots + n_p=n-1$ and that
$$f(v)=(tn_1,tn_2,\ldots,tn_p)$$
for all $v \in V(H)$. In this case, we say that $H$ is a block of type \textbf{(K)}.
\item[(3)] $H=tC_n$ with $t \geq 1$ and $n \geq 5$ odd and there are two indices $k \neq \ell$ from the set $\{1,2,\ldots,p\}$ such that
$$f_i(v)=
\begin{cases}
t & \text{if } i \in \{k,\ell\}, \\
0 & \text{otherwise}
\end{cases}
$$ for all $i \in \{1,2,\ldots,p\}$ and for each $v \in V(H)$. In this case, we say that $H$ is a block of type \textbf{(C)}.
\item[(4)] There are two hard pairs $(H^1,f^1)$ and $(H^2,f^2)$ with $f^1 \in \mathcal{V}_p(H^1)$ and $f^2 \in \mathcal{V}_p(H^2)$ such that $H$ is obtained from $H^1$ and $H^2$ by merging two vertices $v^1 \in V(H_1)$ and $v^ 2 \in V(H_2)$ to a new vertex $v^*$. Furthermore, it holds
$$f(v)=
\begin{cases}
f^1(v) & \text{if } v \in V(H_1) \sm \{v^1\}, \\
f^2(v) & \text{if } v \in V(H_2) \sm \{v^2\}, \\
f^1(v^1) + f^2(v^2) & \text{if } v=v^*
\end{cases}$$
for all $v \in V(H)$. In this case we say that $(H,f)$ is obtained from $(H^1,f^1)$ and $(H^2,f^2)$ by merging $v^1$ and $v^2$ to $v^*$.
\end{itemize}

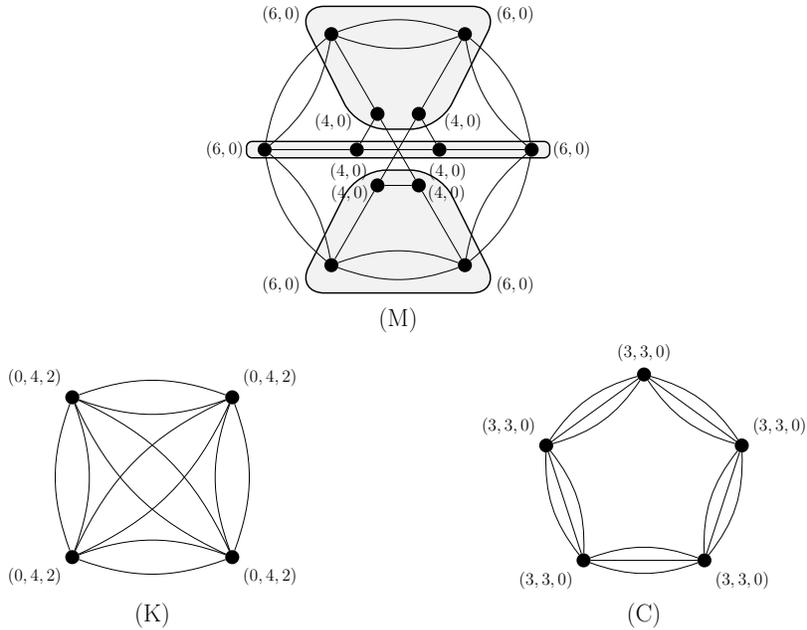
\begin{figure}[thbp]
\centering \resizebox{.81\linewidth}{!}{
\begin{tikzpicture} [node distance=1cm, bend angle=20,
vertex/.style={circle,minimum size=3mm,very thick, draw=black, fill=black, inner sep=0mm}, information text/.style={fill=red!10,inner sep=1ex, font=\Large}, help lines/.style={-,color=black}]

\node[draw=none,minimum size=2cm,regular polygon,regular polygon sides=6, xshift=6cm] (b) {};

\node[vertex] (w1) [label={[xshift=.3cm, yshift=-.2cm]right:{$(4,0)$}}] at (b.corner 1) {};
\node[vertex] (w5) [label={[xshift=-.1cm, yshift=-.2cm]right:{$(4,0)$}}] at (b.corner 5) {};
\node[vertex] (w3) [label={[xshift=-.2cm]below:{$(4,0)$}}] at (b.corner 3) {};
\node[vertex] (w6) [label={[xshift=.2cm]below:{$(4,0)$}}] at (b.corner 6) {};
\node[vertex] (w2) [label={[xshift=-.3cm, yshift=-.2cm]left:{$(4,0)$}}] at (b.corner 2) {};
\node[vertex] (w4) [label={[xshift=.1cm, yshift=-.2cm]left:{$(4,0)$}}] at (b.corner 4) {};

\node[draw=none,minimum size=6.5cm,regular polygon,regular polygon sides=6, xshift=6cm] (c) {};

\node[vertex] (u1) [label={[xshift=.5cm]north east:{$(6,0)$}}] at (c.corner 1) {};
\node[vertex] (u2) [label={[xshift=-.5cm]north west:{$(6,0)$}}]at (c.corner 2) {};
\node[vertex] (u3) [label={[xshift=-.2cm]left:{$(6,0)$}}]at (c.corner 3) {};
\node[vertex] (u4) [label={[xshift=-.5cm]south west:{$(6,0)$}}]at (c.corner 4) {};
\node[vertex] (u5) [label={[xshift=.5cm]south east:{$(6,0)$}}]at (c.corner 5) {};
\node[vertex] (u6) [label={[xshift=.2cm]right:{$(6,0)$}}]at (c.corner 6) {};

\begin{pgfonlayer}{background}
\filldraw[fill=black!50!, fill opacity=.1,rounded corners=20pt, line width=1pt, draw=black] (3.5,3.5) -- (8.5,3.5)--(7,.5)-- (5,.5) -- cycle;

\filldraw[fill=black!50!, fill opacity=.1,rounded corners=20pt, line width=1pt, draw=black] (3.5,-3.5) -- (8.5,-3.5)--(7,-.5)--(5,-.5) -- cycle;

\filldraw[fill=black!50!, fill opacity=.1,  line width=1pt, rounded corners=4pt, draw=black] (2.3,-0.2) rectangle (9.7,0.2);
\end{pgfonlayer}

\node[draw=none, minimum size=5cm, regular polygon, regular polygon sides=5, yshift=-8cm, xshift=12cm](e){};
\node[vertex] (v51) [label={north:$(3,3,0)$}] at (e.corner 1) {};
\node[vertex] (v52) [label={north west:$(3,3,0)$}]at (e.corner 2) {};
\node[vertex] (v53) [label={south west:$(3,3,0)$}]at (e.corner 3) {};
\node[vertex] (v54) [label={south east:$(3,3,0)$}]at (e.corner 4) {};
\node[vertex] (v55) [label={north east:$(3,3,0)$}]at (e.corner 5) {};

\node[draw=none, minimum size=5.5cm, regular polygon, regular polygon sides=4, yshift=-8cm ] (d) {};

\node[vertex] (v41) [label={north east:$(0,4,2)$}] at (d.corner 1) {};
\node[vertex] (v42) [label={north west:$(0,4,2)$}] at (d.corner 2) {};
\node[vertex] (v43) [label={south west:$(0,4,2)$}] at (d.corner 3) {};
\node[vertex] (v44) [label={south east:$(0,4,2)$}] at (d.corner 4) {};
		 
\path[-]
		(w1)
		edge[help lines] (w6)
		edge[help lines] (w4)
		edge[help lines] (u1)
		(w2)
		edge[help lines] (w3)
		edge[help lines] (w5)
		edge[help lines] (u2)
		(w3)
		edge[help lines] (w6)
		edge[help lines] (u3)
		(w4)
		edge[help lines] (w5)
		edge[help lines] (u4)
		(w5)
		edge[help lines] (u5)
		(w6)
		edge[help lines] (u6)
		(u1)
		edge[help lines, bend left] (u2)
		edge[help lines, bend right] (u2)
		(u2)
		edge[help lines, bend left] (u3)
		edge[help lines, bend right] (u3)
		(u3)
		edge[help lines, bend left] (u4)
		edge[help lines, bend right] (u4)
		(u4)
		edge[help lines, bend left](u5)
		edge[help lines, bend right] (u5) 
		(u5)
		edge[help lines, bend left] (u6)
		edge[help lines, bend right] (u6)
		(u6)
		edge[help lines, bend left] (u1)
		edge[help lines, bend right] (u1);

\path [-](v41) edge[help lines, bend left]  (v42)
		   	   edge[help lines, bend right] (v42)
			   edge[help lines, bend left]  (v43)
			   edge[help lines, bend right] (v43)			   
			   edge[help lines,bend left]   (v44)
			   edge[help lines, bend right] (v44)	       
		 (v42) edge[help lines, bend left]  (v43)
			   edge[help lines, bend right] (v43)		 	   
		 	   edge[help lines, bend left]  (v44) 
		 	   edge[help lines, bend right] (v44)    
		 (v43) edge[help lines, bend left]  (v44)
		 	   edge[help lines, bend right] (v44);
		 	   
\path[-]
		 (v51) edge[help lines, bend right] (v52) 
		 (v52) edge[help lines, bend right] (v53) 
		 (v53) edge[help lines, bend right] (v54)
		 (v54) edge[help lines, bend right] (v55) 
		 (v55) edge[help lines, bend right] (v51)
		 (v51) edge[help lines, bend left]  (v52) 
		 (v52) edge[help lines, bend left]  (v53) 
		 (v53) edge[help lines, bend left]  (v54)
		 (v54) edge[help lines, bend left]  (v55) 
		 (v55) edge[help lines, bend left]  (v51)
		 (v51) edge[help lines]  (v52) 
		 (v52) edge[help lines]  (v53) 
		 (v53) edge[help lines]  (v54)
		 (v54) edge[help lines]  (v55) 
		 (v55) edge[help lines]  (v51);
		 
\begin{pgfonlayer}{background}
\node (w41) at (v41) [xshift=1.2cm, yshift=1cm]{};
\node (w42) at (v42) [xshift=-1.2cm, yshift=1cm]{};
\node (w43) at (v43) [xshift=-1.2cm, yshift=-1cm]{};
\node (w44) at (v44) [xshift=1.2cm, yshift=-1cm]{};

\node (w51) at (v51) [yshift=1cm]{};
\node (w52) at (v52) [xshift=-1.2cm]{};
\node (w53) at (v53) [xshift=-1.2cm, yshift=-1cm]{};
\node (w54) at (v54) [xshift=1.2cm, yshift=-1cm]{};
\node (w55) at (v55) [xshift=1.2cm]{};

\node (w11) at (u1) [xshift=1.2cm, yshift=1cm]{};
\node (w12) at (u2) [xshift=-1.2cm, yshift=1cm]{};
\node (w13) at (u3) [xshift=-1.2cm]{};
\node (w14) at (u4) [xshift=-1.2cm, yshift=-1cm]{};
\node (w15) at (u5) [xshift=1.2cm, yshift=-1cm]{};
\node (w16) at (u6) [xshift=1.2cm]{};

\node [rectangle,  fit=(w51)(w52)(w53)(w54)(w55), label={[font=\Large, yshift=.4 cm]below:(C)}]{};
\node [rectangle,  fit=(w11)(w12)(w13)(w14)(w15)(w16), label={[font=\Large, yshift=.4 cm]below:(M)}]{};

\begin{scope}
\node [rectangle, fit=(w51)(w52)(w53)(w54)(w55), xshift=-12cm, label={[font=\Large, yshift=.4 cm]below:(K)}]{};
\end{scope}
\end{pgfonlayer}		
		 
\end{tikzpicture}
}
\caption{Some examples of hard pairs.}
\end{figure}

In the following section, we will show that if $H$ is a hypergraph and $f \in \mathcal{V}_p(H)$ is a function $(p \geq 1)$ the condition $f_1(v) + f_2(v) + \ldots + f_p(v) \geq d_H(v)$ for all $v \in V(H)$ is not sufficient for the existence of an $f$-partition of $H$ if and only if at least one component of $H$ is $f$-hard. Note that $H$ is $f$-partitionable if and only if each component of $H$ is $f$-partitionable. Thus, it is satisfactory to consider only connected hypergraphs. The next result was proven by Borodin, Kostochka and Toft \cite{BorKosToft} for the class of simple graphs. In the next section, we will show how to extend it to hypergraphs.

\begin{theorem}
\label{Theorem:Hauptsatz}
Let $H$ be a connected hypergraph and let $f\in \cV_p(H)$ be a vector function with $p\geq 1$ such that $f_1(v)+f_2(v)+\cdots +f_p(v)\geq d_H(v)$ for all $v\in V(H)$. Then $H$ is not $f$-partitionable if and only if $(H,f)$ is a hard pair.
\end{theorem}

\section{Proof of Theorem \ref{Theorem:Hauptsatz}}
The proof of Theorem~\ref{Theorem:Hauptsatz} is divided into two parts. In the first part, we prove some properties of hard pairs and show that any hard pair is not $f$-partitionable. The proof of the next proposition can be done by induction on the number of blocks of $H$ and is straightforward.

\begin{proposition}\label{prop_block-hardpair}
Let $H$ be a connected hypergraph, and let $f \in \mathcal{V}_p(H)$ be a vector function with $p \geq 1$ such that $H$ is $f$-hard.  Then, for each $B \in \cB(H)$ there is a uniquely determined function $f_B \in \mathcal{V}_p(B)$ such that the following statements hold.
\begin{itemize}
\item[\upshape (a)] $(B,f_B)$ is a hard pair of type {\upshape (M), (K)} or {\upshape(C)}.
\item[\upshape (b)] $f(v) = \sum_{B \in \cB_v(H)}f_B(v)$ for all $v \in V(H)$.
\item[\upshape (c)] $f_B(v)=f(v)$ for all non-separating vertices of $H$ belonging to $B$.
\end{itemize}
\end{proposition}

 Note that the above proposition clearly implies that $f_B(v) \leq f(v)$ holds componentwise. The next proposition shows that $f$-hard hypergraphs are not $f$-partitionable.

\begin{proposition} \label{prop_f-hard}
Let $H$ be a connected hypergraph, and let $f \in \mathcal{V}_p(H)$ be a vector function with $p \geq 1$. If $H$ is $f$-hard, then the following statements hold.
\begin{itemize}
\item[\upshape(a)] $f_1(v) + f_2(v) + \ldots + f_p(v) = d_H(v)$ for all $v \in V(H)$.
\item[\upshape(b)] If $u \neq u'$ are two non-separating vertices contained in the same block of $H$, then either $f(u)=f(u')$ or $f_i(u)=f_i(u')=0$ for all but one index $i \in \{1,2,\ldots,p\}$.
\item[\upshape(c)] $H$ is not $f$-partitionable.
\end{itemize}
\end{proposition}

\begin{proof}
Statements (a) and (b) are simple consequences of Proposition~\ref{prop_block-hardpair}. The proof of (c) is by reductio ad absurdum. To this end, choose $(H,f)$ such that
\begin{itemize}
\item[(1)] $H$ is $f$-hard,
\item[(2)] there is an $f$-partition $(H_1,H_2,\ldots,H_p)$ of $H$, and
\item[(3)] $|H|$ is minimum with respect to (1) and (2).
\end{itemize}
Note that the empty hypergraph is the only hypergraph that is strictly $0$-degenerate; thus, if $f_i \equiv 0$ for some $i$, then $H_i=\ems$ must hold. As a consequence, if $(H,f)$ is of type (M), there is an index $j$ such that $H_i=\ems$ for all $i \in \{1,2,\ldots,p\} \setminus \{j\}$ and $f_j(v)=d_H(v)$ for all $v \in V(H)$. Therefore, $H_j$ is not strictly $f_j$-degenerate, contradicting (2).

If $(H,f)$ is of type (K), then $H=tK_n$ for some $t \geq 1, n \geq 3$ and there are integers $n_1,n_2,\ldots,n_p$ such that $n_1+n_2+\ldots+n_p=n-1$ and $f(v)=(tn_1,tn_2,\ldots,tn_p)$ for all $v \in V(H)$. Thus, $H_i$ is a $tK_{m_i}$ for some $m_i \geq 0$ for all $i \in \{1,2,\ldots,p\}$. Since $H_i$ is strictly $f_i$-degenerate, it holds $|H_i| \leq n_i$ for all $i \in \{1,2,\ldots,p\}$. Consequently, we obtain $$|H| = |H_1| + |H_2| + \ldots + |H_p| \leq n-1,$$ which is impossible.

If $(H,f)$ is of type (C), then $H=tC_n$ for some $t \geq 1$ and $n \geq 5$ odd, and there are two indices $k \neq \ell$ from the set $\{1,2,\ldots,p\}$ such that $f_i(v)=t$ for $i \in \{k,\ell \}$ and $f_i(v)=0$, otherwise. Then, $(H_k,H_\ell)$ is a 2-partition of $H$ and $f_k(v)=f_\ell(v)=t$ for all $v \in V(H)$. Since $n$ is odd, one of the parts, say $H_k$, contains two adjacent vertices that are joined by $t$ parallel edges. Therefore, $H_k$ is not strictly $f_k$-degenerate, a contradiction.

It remains to consider the case that $(H,f)$ is obtained from two hards pairs $(H^1,f^1)$ and $(H^2,f^2)$ by merging $v^1$ and $v^2$ to $v^*$. By (3), $H^j$ is not $f^j$-partitionable for $j \in \{1,2\}$. Let $H_i^j=H^j\cap H_i$ for $i \in \{1,2,\ldots,p\}$ and $j \in \{1,2\}$. By symmetry, we may assume $v^* \in V(H_1)$. Since $(H_1,H_2,\ldots,H_p)$ is an $f$-partition of $H$, it follows that $H_i^1$ is strictly $f_i^1$-degenerate and $H_i^2$ is strictly $f_i^2$-degenerate for all $i \in \{2,3,\ldots,p\}$. As a consequence, for $j \in \{1,2\}$, the hypergraph $H_1^j$ is not strictly $f_1^j$-degenerate and, thus, there is a non-empty subhypergraph $G^j \subseteq H_1^j$ such that $d_{G^j}(v) \geq f_1^j(v)$ for all $v \in V(G^j)$. Nevertheless, this implicates that $G=G^1\cup G^2$ is a non-empty subhypergraph of $H_1$ such that $d_G(v) \geq f_1(v)$ for all $v \in V(G)$, a contradiction. This completes the proof.
\end{proof}

Thus, the 'if'-direction is proved. For the remaining part, we will need the following notation. We say that $(H,f)$ is a \textbf{non-partionable pair of dimension} $p$ if $H$ is a connected hypergraph, $f \in \mathcal{V}_p(H)$ is a vector function satisfying $$f_1(v) + f_2(v) + \ldots + f_p(v) \geq d_H(v)$$ for all $v \in V(H)$, and $H$ is not $f$-partitionable. The next two propositions describe characteristics of non-partitionable pairs.

\begin{proposition}\label{prop_reduction}
Let $(H,f)$ be a non-partitionable pair of dimension $p$, let $z$ be a non-separating vertex of $H$, and let $j \in \{1,2,\ldots,p\}$ such that $f_j(z) \neq 0$. For the hypergraph $H'=H \div z$, define $f' \in \mathcal{V}_p(H')$ to be the vector function satisfying
$$f_i'(v)=\begin{cases}
\max\{0, f_j(v) - \mu_H(z,v)\} & \text{if } i=j, \\
f_i(v) & \text{otherwise}
\end{cases}$$
for all $v \in V(H')$ and $i \in \{1,2,\ldots,p\}$. Then, $(H',f')$ is a non-partition\-able pair of dimension $p$, and in what follows, we write $(H',f')=(H,f)/(z,j)$.
\end{proposition}

\begin{proof}
By symmetry, we may assume $j=1$. Then, $f_1(z) \geq 1$ and $H$ is not $f$-partitionable. Thus, $|H| \geq 2$ holds and $H'=H \div z$ is connected. Assume that $H'$ admits an $f'$-partition $(H_1,H_2, \ldots, H_p)$. To arrive at a contradiction, let $H_1^*=H[V(H_1) \cup \{z\}]$. We show that $H_1^*$ is strictly  $f_1$-degenerate. To this end, choose a non-empty subhypergraph $G^* \subseteq H_1^*$.  Then, $H_1=H_1^* \div z$ and $G=G^* \div z$ is a subhypergraph of $H_1$. As $H_1$ is strictly $f_1'$-degenerate, $G$ is strictly $f_1'$-degenerate, too. If $G$ is non-empty, this implies that there is a vertex $v$ satisfying $ d_G(v) < f_1'(v)$. But then, $f_1'(v) > 0$ and, by using \eqref{Equ:degree(uinH-v)}, we obtain
$$d_{G^*}(v) = d_G(v) + \mu_H(v,z) < f_1'(v) + \mu_H(v,z) = f_1(v),$$
and we are done. If $G$ is empty, then $V(G^*)=\{z\}$ and $d_{G^*}(z)=0<f_1(z)$. Hence, $H_1^*$ is strictly $f_1$-degenerate. As $H_j^*=H[V(H_j)]$ is a subhypergraph of $H_j$ for $j \in \{2,3,\ldots,p\}$, $H_j^*$ is strictly-$f_j$-degenerate and, thus, the sequence $(H_1^*,H_2^*,\ldots,H_p^*)$ is an $f$-partition of $H$, which is impossible.
\end{proof}

By applying the above introduced reduction method, we obtain the following statements.

\begin{proposition}\label{prop_degree-condition}
Let $(H,f)$ be a non-partitionable pair of dimension $p \geq 1$. Then, the following statements hold.
\begin{itemize}
\item[\upshape (a)] $f_1(v) + f_2(v) + \ldots + f_p(v) = d_H(v)$ for all $v \in V(H)$.
\item[\upshape (b)] If $z$ is a non-separating vertex of $H$ satisfying $f_j(z) \neq 0 $ for some $j \in \{1,2,\ldots,p\}$, then $f_j(v) \geq \mu_H(z,v)$ holds for all $v \in V(H) \setminus\{z\}$.
\item[\upshape (c)] If $|H| \geq 2$ and if $u$ is an arbitrary vertex of $H$, then $H-u$ admits an $f$-partition. Furthermore, for any $f$-partition $(H_1,H_2,\ldots,H_p)$ it holds $f_i(u)=d_{H_i+u}(u)$ for all $i \in \{1,2,\ldots,p\}$ and $E_H(u) = E_{H_1+u}(u) \cup E_{H_2 + u}(u) \cup \ldots \cup E_{H_p + u}(u)$.
\end{itemize}
\end{proposition}

\begin{proof}
The proof of statement (a) is by induction on the order $n$ of $H$. For $n=1$, the statement is evident. Let $n \geq 2$ and let $v$ be an arbitrary vertex. Since $H$ is connected, there is a non-separating vertex $z \neq v$ in $H$. Since $$f_1(z) + f_2(z) + \ldots + f_p(z) \geq d_H(z) \geq 1,$$ it holds $f_j(z) \geq 1$ for some $j \in \{1,2,\ldots,p\}$. By Proposition~\ref{prop_reduction}, the pair $(H',f')=(H,f)/(z,j)$ is non-partitionable, $f_i'(v)=f_i(v)$ for all $i \neq j$ from the set  $\{1,2,\ldots,p\}$ and $f_j'(v)=\max\{0,f_j(v) - \mu_H(v,z)\}$. From the induction hypothesis it follows that
\begin{align*}
f_1'(v) + f_2'(v) + \ldots + f_p'(v)=d_{H'}(v).
\end{align*}
Since $f_1(v) + f_2(v) + \ldots + f_p(v) \geq d_H(v)$, this leads to
\begin{align*}
d_H(v) & \leq f_1(v) + f_2(v) + \ldots + f_p(v)\\
 & \leq f_1'(v) + f_2'(v) + \ldots + f_p'(v) + \mu_H(v,z)\\
 & = d_{H'}(v) + \mu_H(v,z) = d_H(v)
\end{align*}
(see \eqref{Equ:degree(uinH-v)}), and the proof of (a) is complete.

The proof of (b) is by contradiction. Assume that there exist a non-separating vertex $z$ of $H$ and a vertex $v \neq z$ such that $f_j(z) \neq 0$ and $f_j(v) < \mu_H(v,z)$ for some $j \in \{1,2,\ldots,p\}$. By symmetry, we may assume $j=1$. Then, $(H',f')=(H,f)/(z,1)$ is a non-partitionable pair such that $$f_1(v) - \mu_H(z,v)< 0 = f_1'(v)$$ and $f_i(v)=f_i'(v)$ for all $i \in \{2,3,\ldots,p\}$. Using \eqref{Equ:degree(uinH-v)} and applying (a) to $(H',f')$ as well as $(H,f)$ leads to
\begin{align*}
d_H(v) - \mu(z,v)
& = d_{H'}(v) =  f'_1(v) + f'_2(v) + \ldots + f'_p(v)\\
& >  f_1(v) - \mu_H(z,v) + f_2(v) + \ldots + f_p(v)\\
& =  d_H(v) - \mu_H(z,v),
\end{align*}
which is impossible.

In order to prove (c), let $u$ be an arbitrary vertex of $H$ and let $H'=H-u$. Since $H$ is connected, each component $G$ of $H'$ contains a vertex $u'$, which is a neighbor of $u$ in $H$ and, so, $f_1(u') + f_2(u') + \ldots + f_p(u') \geq d_H(u') > d_G(u')$. Applying (a) to $(G,f)$, this implies that $G$ is $f$-partitionable and, thus, $H'$ is $f$-partitionable. Hence, there is an $f$-partition $(H_1,H_2,\ldots,H_p)$ of $H'$. Since $H$ is not $f$-partitionable, we conclude that $H_i'=H_i + u$ is not strictly $f_i$-degenerate for each $i \in \{1,2,\ldots,p\}$. Hence, there is a non-empty subhypergraph $G_i$ of $H_i'$ such that $d_{G_i}(v) \geq f_i(v)$ for all $v \in V(G)$. Since $H_i$ is strictly $f_i$-degenerate, $u \in V(G_i)$ for all $i \in \{1,2,\ldots,p\}$. Due to the fact that $f_1(u) + f_2(u) + \ldots + f_p(u) =d_H(u)$ (by (a)) and $$d_H(u) \geq d_{H_1'}(u) + d_{H_2'}(u) + \ldots + d_{H_p'}(u),$$ it follows that
\begin{align*}
f_1(u) + f_2(u) + \ldots + f_p(u) & = d_H(u)\\
& \geq d_{H_1'}(u) + d_{H_2'}(u) + \ldots + d_{H_p'}(u)\\
& \geq d_{G_1}(u) + d_{G_2}(u) + \ldots + d_{G_p}(u)\\
& \geq f_1(u) + f_2(u) + \ldots + f_p(u),\\
\end{align*}
which leads to $f_i(u)=d_{H_i'}(u)$ for all $i \in \{1,2,\ldots,p\}$. Furthermore, it follows $$d_H(u) = d_{H_1'}(u) + d_{H_2'}(u) + \ldots + d_{H_p'}(u),$$ which clearly implies the last part of the statement.
\end{proof}

Now we are able to prove the remaining part of Theorem~\ref{Theorem:Hauptsatz}.

\begin{theorem}
If $(H,f)$ is a non-partitionable pair of dimension $p \geq 1$, then $H$ is $f$-hard.
\end{theorem}

\begin{proof}
The proof is by reductio ad absurdum. So let $(H,f)$ be a smallest counterexample, that is,
\begin{itemize}
\item[\upshape (1)] $(H,f)$ is a non-partitionable pair of dimension $p \geq 1$,
\item[\upshape (2)] $(H,f)$ is no hard pair, and
\item[\upshape (3)] $|H|$ is minimum subject to (1) and (2).
\end{itemize}
From Proposition~\ref{prop_degree-condition}(a) it then follows that
\begin{align} \label{eq_prop2.5}
f_1(v) + f_2(v) + \ldots + f_p(v) = d_H(v)
\end{align}
for all $v \in V(H)$. Furthermore, $|H|\geq 2$, for otherwise, $(H,f)$ would be a hard pair of type (M), contradicting (2). To arrive at a contradiction, we shall establish seven claims analyzing the structure of the pair $(H,f)$.

\begin{claim}\label{claim_H is block}
$H$ is a block, that is, $H$ has no separating vertex.
\end{claim}

\begin{proof2}
Suppose, to the contrary, that $H$ has a separating vertex $v^*$. Then, $H$ is the union of two connected induced subhypergraphs $H^1$ and $H^2$ with $V(H^1) \cap V(H^2) = \{v^*\}$ and $|H^j| < |H|$ for $j \in \{1,2\}$. By Proposition~\ref{prop_degree-condition}(c), $H-v^*$ admits an $f$-partition $(H_1,H_2, \ldots, H_p)$ satisfying $f_i(v^*)=d_{H_i+v^*}(v^*)$ for all $i \in \{1,2,\ldots,p\}$ and $E_H(v^*) = E_{H_1+v^*}(v^*) \cup E_{H_2 + v^*}(v^*) \cup \ldots \cup E_{H_p + v^*}(v^*)$. For $i \in \{1,2,\ldots,p\}$, we define $H_i^1=H_i \cap H^1$ and $H_i^2=H_i \cap H^2$. Then, $H_i=H_i^1 \cup H_i^2$ and
\begin{align}\label{eq_H_i+v^*}
f_i(v^*)=d_{H_i+v^*}(v^*)=d_{H_i^1+v^*}(v^*)+d_{H_i^2+v^*}(v^*)
\end{align}
for all $i \in  \{1,2,\ldots,p\}$. For $j \in \{1,2\}$ let $f^j \in \mathcal{V}_p(H^j)$ be the function satisfying
$$f_i^j(v)=\begin{cases}
f_i(v) & \text{if } v \in V(H^j-v^*),\\
d_{H_i^j+v^*}(v^*) & \text{if } v=v^*
\end{cases}$$
for all $v \in V(H^j)$ and all $i \in \{1,2,\ldots,p\}$. By \eqref{eq_prop2.5} and \eqref{eq_H_i+v^*} together with Proposition~\ref{prop_degree-condition}(c), we conclude that $f_1^j(v) + f_2^j(v) + \ldots + f_p^j(v) = d_{H^j}(v)$ for each $j \in \{1,2\}$ and $v \in V(H^j)$.
If $H^j$ is not $f^j$-partitionable for $j \in \{1,2\}$, then, as $(H^j,f^j)$ satisfies (1) and since $|H^j| < |H|$, it follows from (3) that $H^j$ is $f^j$-hard. Therefore, $(H,f)$ is obtained from two hard pairs by merging two vertices, and so $H$ is $f$-hard. Otherwise, by symmetry, we may assume that $H^1$ admits an $f^1$-partition $(H_1',H_2',\ldots,H_p')$ and that $v^* \in V(H_1')$. Consider the $p$-partition $(G_1,G_2,\ldots,G_p)$ of $H$, whereby $G_1=H_1' \cup (H_1^2+v^*)$ and $G_i = H_i' \cup H_i^2$ for $i \in \{2,3,\ldots,p\}$. By construction, $G_i$ is strictly $f_i$-degenerate for $i \in \{2,3,\ldots,p\}$. We claim that $G_1$ is strictly $f_1$-degenerate. In order to prove this, let $G$ be a non-empty subhypergraph of $G_1$. If $G \subseteq H_1^2$, then $d_G(v) < f(v)$ holds for some vertex $v \in V(G)$, since $H_1^2$ is strictly $f_1$-degenerate. Otherwise, $G'=G \cap H_1'$ is a non-empty subhypergraph of $H_1'$ and, since $H_1'$ is strictly $f_1^1$-degenerate, there is a vertex $v \in V(G')$ such that $d_{G'}(v) < f_1^1(v)$. If $v \neq v^*$, then $d_G(v) = d_{G'}(v) < f_1^1(v) = f_1(v)$ and we are done. Else, $v=v^*$ and it follows from \eqref{eq_H_i+v^*} and from the definition of $f_1^j$ that
$$ d_G(v^*) \leq d_{G'}(v^*) + d_{H_1^2+v^*}(v^*) < f_1^1(v^*) + f_1^2(v^*) = f_1(v^*).$$
This shows that $G_1$ is strictly $f_1$-degenerate and, hence, $H$ is $f$-partitionable, contradicting the premise. Thus, the first case is complete.
\end{proof2}

\begin{claim}\label{claim_multiplicity}
If there exists a vertex $z \in V(H)$ and an index $j \in \{1,2,\ldots,p\}$ such that $f_j(z) \neq 0$, then $(H',f')=(H,f)/(z,j)$ is  a non-partitionable pair and the following statements hold.
\begin{itemize}
\item[\upshape (a)] $(H',f')$ is a hard pair.
\item[\upshape (b)] $f_j(v) \geq \mu_H(v,z)$ for all $v \in V(H)\setminus \{z\}$.
\end{itemize}
\end{claim}

\begin{proof2}
Since $H$ is a block (by Claim~\ref{claim_H is block}) and $|H|\geq 2$, $z$ is a non-separating vertex of $H$ and $H'=H\div z\not=\ems$. Since $(H,f)$ is a non-partionable pair (by (1)), $(H',f')$ is a non-partionable pair, too (by Proposition~\ref{prop_reduction}). From (3) it then follows that $(H',f')$ is a hard pair. Statement (b) is a consequence of Proposition~\ref{prop_degree-condition}(b).
\end{proof2}

Now let $z\in V(H)$ be an arbitrary vertex. Since $|H| \geq 2$ and since $H$ is connected, there is an index $j \in \{1,2,\ldots,p\}$ with $f_j(z)\neq 0$ (by \eqref{eq_prop2.5}). By symmetry, we may assume $j=1$. Then, $(H',f')=(H,f)/(z,1)$ is a hard pair (by Claim~\ref{claim_multiplicity}(a)). Furthermore, $f_1(v) \geq \mu_H(v,z)$ (by Claim~\ref{claim_multiplicity}(b)), and so \begin{align} \label{eq_degreefirstcase}
f'(v)=(f_1(v) - \mu_H(z,v),f_2(v),\ldots,f_p(v))
\end{align}
for all $v \in V(H')$.

\begin{claim}\label{claim_mono}
The hard pair $(H',f')$ is not of type {\upshape (M)}.
\end{claim}
\begin{proof2}
Assume, by contrary, that $(H',f')$ is of type (M). If $n=2$ this implies that $(H,f)$ is of type (M) since otherwise $(H,f)$ would clearly admit an $f$-partition. But then $(H,f)$ is a hard pair, contradicting (2). Now let $n \geq 3$. Since $(H',f')$ is of type (M), there is an index $j \in \{1,2,\ldots,p\}$ such that $f_j'(v)=d_{H'}(v)$ and $f_i'(v)=0$ for all $i \in \{1,2,\ldots, p\} \setminus \{j\}$ and for all $v \in V(H')$.
\medskip

\case{A}{The vertex $z$ is contained in an hyperedge.} Then, $H-z$ is a proper subhypergraph of $H'=H \div z$. Since each vertex $v$ of $H'$ satisfies $d_{H'}(v)=f_j'(v)$, Proposition~\ref{prop_reg-f-part} implies that $H-z$ is strictly $f_j'$-degenerate and therefore strictly $f_j$-degenerate. If $j \neq 1$, then setting $H_1=H[\{z\}],~H_j=H-z$, and $H_i=\ems$ for $i \in \{2,3,\ldots,p\} \setminus \{j\}$ gives us an $f$-partition of $H$, which is impossible. Thus, $j=1$. Moreover, by a similar argumentation it must hold $f(z)=(d_H(z),0,\ldots,0)$ and, thus, $(H,f)$ is a hard pair of type (M), contradicting (2).
\medskip

\case{B}{The vertex $z$ is contained only in ordinary edges.} Since $H$ is a block, this implies that the set $$N = \{v \in V(H') ~ | ~ vz \in E(H)\}$$ is non-empty. If $j=1$, then \eqref{eq_degreefirstcase} leads to $f_2(v)=f_3(v)=\ldots=f_p(v)=0$ for all $v \in V(H')$. By Claim~\ref{claim_multiplicity}(b), it then follows $f_2(z)=f_3(z)=\ldots=f_p(z)=0$ and, thus, $(H,f)$ is of type (M), a contradiction to (2).

It remains to consider the case that $j \neq 1$, say $j=2$ (by symmetry). Then, $f_3(v)=f_4(v)=\ldots = f_p(v) = 0$ and $f_2(v)=f_2'(v)=d_{H'}(v)>0$ for all $v \in V(H')$. Since $N$ is non-empty, Claim~\ref{claim_multiplicity}(b) (applied on a vertex from $N$) implies that $f_2(z) > 0$. Then, $(H',f'')=(H,f)/(z,2)$ is a hard pair (by Claim~\ref{claim_multiplicity}(a)), too, and it holds $$f''(v)=(f_1(v), d_{H'}(v) - \mu_H(v,z),0,0,\ldots,0)$$ for all $v \in V(H')$. Assume that there is a vertex $u \in V(H') \setminus N$. Then, $\mu_H(u,z)=0$ and $f''(u)=f'(u)=(0,d_{H'}(u),0,0,\ldots,0)$. Since $(H',f'')$ is a hard pair, we conclude that $(H',f'')$ is a hard pair of type (M) with $f''(v)=(0,d_{H'}(v),0,0,\ldots,0)$ for all $v \in V(H')$. However, since $f_1(z) > 0$, Claim~\ref{claim_multiplicity}(b) leads to $f_1(v) > 0$ for all $v \in N$ and, thus, $f_1''(v)=f_1(v) > 0$ for all $v \in N$, a contradiction. As a consequence, $N=V(H')$. Then, $f_1(z) > 0$ and $f_2(v)=f_2'(v)=d_{H'}(v)>0$ for all $v \in V(H')$. By Claim~\ref{claim_multiplicity}(b), this leads to $f_1$ and $f_2$ being nowhere-zero in $V(H)$. Let $v \in V(H')$ be an arbitrary vertex. Then, $H_2=H-z-v$ is a proper subhypergraph of $H'$ and, therefore, strictly $f_2'$-degenerate (by Proposition~\ref{prop_reg-f-part}). If $H_1=H[\{v,z\}]$ is strictly $f_1$-degenerate, then $(H_1,H_2,\ems,\ems,\ldots,\ems)$ is an $f$-partition of $H$, which is impossible. Thus, $H_1=H[\{v,z\}]$ is not strictly $f_1$-degenerate and, since $f_1(z), f_1(v) \geq \mu_H(v,z) \geq 1$ (by Claim~\ref{claim_multiplicity}(b)), this leads to $f_1(z)=f_1(v)=\mu_H(v,z)$. Since $v$ was chosen arbitrarily, this implies that there is an integer $m \geq 1$ such that $m=\mu_H(v,z)=f_1(v)=f_1(z)$ for all $v \in V(H')=N$. Since $n \geq 3$, $N$ contains at least two vertices. We choose two different vertices from $N$, say $u$ and $v$ and show that $\mu_H(u,v)=m$. Let $H_2=H-u-v$. We claim that $H_2$ is strictly $f_2$-degenerate. To this end, let $G$ be a non-empty subhypergraph of $H_2$. If $z$ is contained in $G$, then $$d_G(z) \leq d_{H_2}(z) = d_H(z) - 2m < d_H(z) - m = f_2(z),$$ (as $f_1(z)=m$ and as $f_3(z)=f_4(z)=\ldots=f_p(z)=0)$ and we are done. If $z$ is not contained in $G$, then $G$ is a proper subhypergraph of $H'=H \div z = H-z$. As $d_{H'}(v)=f_2'(v)$ for all $v \in V(H')$, it follows from Proposition~\ref{prop_reg-f-part} that $G$ is strictly $f_2'$-degenerate and, therefore, strictly $f_2$-degenerate. Consequently, $H_2$ is strictly $f_2$-degenerate. If $H_1=H[\{u,v\}]$ is strictly $f_1$-degenerate, then $(H_1,H_2,\ems,\ems,\ldots,\ems)$ is an $f$-partition of $H$, which is impossible. Thus, $H_1$ is not strictly $f_1$-degenerate and, as $m=f_1(u)=f_1(v)\geq \mu_H(u,v)$ (by Claim~\ref{claim_multiplicity}(b)), it must hold $\mu_H(u,v)=m$.

Since $u$ and $v$ were chosen arbitrarily from $N$, we conclude that $\mu_H(u,v)=m$ for all $u \neq v$ from $V(H)$. As a consequence, we obtain $f_2(v) \geq m(n-2)$ for all $v \in V(H')$. Assume that $H$ contains a hyperedge $e$. Then, since $z$ is contained only in ordinary edges, it must hold $e \in E(H')$. Moreover, regarding $(H',f'')=(H,f)/(z,2)$, it follows that $(H',f'')$ is a hard pair of type (M) (by Claim~\ref{claim_multiplicity}(a) and as $H'$ contains $e$). Furthermore, as $f_1(v)=f_1''(v) > 0$ for all $v \in V(H')$, it must hold $f_2''(v)=0$ for all $v \in V(H')$ and, therefore, $f_2(v) = m(n-2)$ for all $v \in V(H')$ and $n=3$. However, this leads to $|H'|=n-1=2$ and, thus, $H'$ cannot contain any hyperedges, a contradiction. Hence, $H$ does not contain any hyperedges and, therefore, $H$ is a $mK_n$ and it holds $f_2(v)=m(n-2)$ for all $v \in V(H')$. As $d_H(z)=m(n-1)=f_1(z) + f_2(z)= m + f_2(z)$ (by \eqref{eq_prop2.5}), we conclude that $f_2(z) = m(n-2)$ and $(H,f)$ is a hard pair of type (K), contradicting (2). This completes the proof.
\end{proof2}

\begin{claim} \label{claim_complete}
The hard pair $(H',f')$ is not of type {\upshape (K)}.
\end{claim}

\begin{proof2}
Assume, by contrary, that $(H',f')=(H,f)/(z,1)$ is of type (K). Then, it holds $H'=tK_{n-1}$ for some $t \geq 1$ and $n \geq 4$ and there are integers $n_1,n_2,\ldots,n_p$ with at least two $n_i$ different from zero such that $n_1+n_2+\ldots+n_p=n-2$ and that $f'(v)=(tn_1,tn_2,\ldots,tn_p)$ for all $v \in V(H')$. By symmetry, we may assume $n_2 > 0$ and, thus, $f_2'(v) > 0$ for all $v \in V(H')$. We distinguish between two cases.

\medskip
\case{A}{$E_H(z)$ contains an ordinary edge.} Then, the set $$N=\{ v \in V(H') ~ | ~ vz \in E(H)\}$$ is non-empty. Since $f_2'(v) > 0$ for all $v \in N$, Claim~\ref{claim_multiplicity}(b) implies that $f_2(z) > 0$. Let $v \in N$ and let $\mu_H(v,z)=m$. Then, $f(v)=(tn_1+m,tn_2,tn_3,\ldots,tn_p)$ (by \eqref{eq_degreefirstcase}). Since $f_2(z) > 0$, Claim~\ref{claim_multiplicity}(a) implies that $(H',f'')=(H,f)/(z,2)$ is also a hard pair, which can be only of type (M) or (K). Since $H'$ is regular, this implies that $f''$ is constant. Furthermore, $$f''(v)=(tn_1+m, tn_2-m, tn_3, \ldots, tn_p) \neq f'(v)$$ and $$f''(w)=(tn_1 + \mu_H(w,z), tn_2 - \mu_H(w,z), tn_3, \ldots, tn_p)$$ for all $w \in V(H') \setminus \{v\}$ (by Claim~\ref{claim_multiplicity}(b)). Consequently, $\mu_H(z,u)=m \geq 1$ for all $u \in V(H')$. Since $(H',f'')$ is of type (M) or (K) and since all vertices of $H'$ have degree $t(n-2)$ in $H'$, we furthermore conclude that $t|m$ and so $m \geq t$. Finally, we obtain that $f_1$ as well as $f_2$ are nowhere-zero in $V(H)$.

Next we claim that $H \div v$ is a block for all $v \in V(H)$. Otherwise, there would exist a vertex $v \in V(H)$ different from $z$ such that $H \div v$ is not a block. Since $z$ is joined to all other vertices by ordinary edges, $z$ would be the only possible separating vertex of $H \div v$. However, as $(H\div z) \div v = (H \div v) \div z$, the hypergraph $(H \div v) \div z$ is complete and therefore connected, a contradiction.

Assume that there is a hyperedge $e \in E(H)$. Then, $z \in i_H(e)$ and $|i_H(e)|=3$ since $H'$ does not contain any hyperedges. Since $n \geq 4$,  there is a vertex $x \in V(H) \setminus i_H(e)$. As $H \div x$ is a block containing the hyperedge $e$, the hard pair $(H,f)/(x,j)$ must be of type (M) for $j \in \{1,2\}$. Since $\mu_H(x,z)=m$ and since $f_1(z), f_2(z) \geq m$ (by Claim~\ref{claim_multiplicity}(b)), this implies that $f_1(z)=f_2(z)=m$ and $f_3(z)=f_4(z)=\ldots=f_p(z)=0$. As a consequence, $f_1(z) + f_2(z) + \ldots + f_p(z) = 2m < 3m \leq d_H(z)$, a contradiction.

Hence, there are no hyperedges in $H$. Then, $H'=H \div z = H-z$ and so $\mu_H(u,v)=t$ for all $u \neq v$ from $V(H) \setminus \{z\}$. Moreover, as $\mu_H(v,z)=m$ for all $v \in V(H) \setminus \{z\}$ and since $f_1$ and $f_2$ are nowhere-zero, it holds $f_1(v) \geq m$ and $f_2(v) \geq m$ (by Claim\ref{claim_multiplicity}(b)). We show that $t=m$ and so $H=tK_n$. Otherwise, $t < m$. As $n \geq 4$, $(H^*,f^*)=(H,f)/(x,1)$ must be of type (M) for any $x \in V(H) \setminus \{z\}$.  However, since $t < m$ it holds $f_1^*(v) =f_1(v) - t> 0$ and $f_2^*(v) = f_2(v) > 0$ for all $v \in V(H^*) \setminus \{z\} \neq \ems$, a contradiction. Thus, $m=t$ and so $H=tK_n$.

To conclude the case, we show that $(H,f)$ is of type (K), giving a contradiction to statement (2). To this end, choose two distinct vertices $u$ and $v$ in $H$. By Proposition~\ref{prop_degree-condition}(c), $H-u$ admits an $f$-partition $(H_1,H_2,\ldots,H_p)$ and $f_i(u)=d_{H_i + u} (u) = t |H_i|$ for every $i \in \{1,2,\ldots,p\}$. By symmetry, we may assume $v \in V(H_1)$. Due to the fact that $H_1$ is strictly $f_1$-degenerate and since $f_1(u) = d_{H_1 + u} (u) > d_{(H_1-v)+u} (u)$, the hypergraph $H_1'=(H_1 - v) + u$ is also strictly $f_1$-degenerate. Thus, $(H_1',H_2,\ldots,H_p)$ is an $f$-partition of $H-v$ satisfying $|H_1'|=|H_1|$. As a consequence, $f_i(v)=t|H_i|$ for every $i \in \{1,2,\ldots,p\}$ (by Proposition~\ref{prop_degree-condition}(c)). In conclusion, $f_i(u)=f_i(v)=t|H_i|$ for each $i \in \{1,2,\ldots,p\}$ and, since $f_1$ and $f_2$ are nowhere-zero, at least two $|H_i|$ are non-empty. Therefore, $(H,f)$ is of type (K), contradicting (2).

\medskip
\case{B}{$E_H(z)$ contains only hyperedges.} This implies that $f(v)=f'(v)=(tn_1,tn_2, \ldots, tn_p)$ for all $v \in V(H')$. First assume that there is a vertex $v \in V(H')$ such that $H \div v$ has a separating vertex. Since
$(H \div v) \div z = (H \div z) \div v$ is complete, $z$ is a non-separating vertex of $\tilde{H}=H \div v$. Let $B$ be the block of $\tilde{H}$ containing $z$. Due to the fact that any two distinct vertices of $V(H) \setminus \{v,z\}$ are either contained in an ordinary edge of $H$ or in an hyperedge of $H$ together with $z$, they all are contained in the same block $B'$ of $\tilde{H}$ and $B'$ is a $tK_{n-2}$. Since $\tilde{H}$ has at least two blocks, this implies that $B$ and $B'$ are the only blocks of $\tilde{H}$ and that there is exactly one separating vertex $u$ in $\tilde{H}$. Moreover, we conclude that there is a hyperedge $e$ in $H$ with $i_H(e)=\{u,v,z\}$. Let $x$ be a non-separating vertex of $B'$.
Then, $(H'',f'')=(H,f)/(x,2)$ is a hard pair (since $f_2(x)=f_2'(x) > 0$ and by Claim~\ref{claim_multiplicity}(a)). As $B'$ is a $tK_{n-2}$ and as $v$ is joined to all vertices from $V(B') \setminus \{u\}$ by ordinary edges (since $H \div z$ is a $tK_{n-1}$), we conclude that $H''$ is a block, which contains the hyperedge $e$. Thus, $(H'',f'')$ is of type (M) and there is an index $i \in \{1,2,\ldots,p\}$ such that $f_i''(w)=d_{H''}(w) > 0$ and $f_k''(w)=0$ for all $k \in \{1,2,\ldots,p\} \setminus \{i\}$ and for all $w \in V(H'')$. In particular, since $f_1(z)=f_1''(z) > 0$, it holds $i=1$. Thus, $f_1(w)=f_1''(w)=d_{H''}(w) > 0$ and $f_k''(w)=0$ for all $k \in \{2,3,\ldots,p\}$ and for all $w \in V(H'')$. However, this also implies that $f_1(x) > 0$ (since $f_1=f_1'$ is constant in $V(H')$). By Claim~\ref{claim_multiplicity}(a),  $(H'',f^*)=(H,f)/(x,1)$ again is a hard pair with $f_1(w) = d_{H''} (w)$ and $f_k(w) = 0$ for all $k \in \{2,3,\ldots,p\}$ and for all $w \in V(H')$ (as $f_1(z)=f_1^*(z)>0$) . However, in this case we obtain $f^*_2(v)=f_2(v) > 0$, a contradiction.

It remains to consider the case that $H \div v$ is a block for all $v \in V(H)$. If $f(z)=(d_H(z),0,0,\ldots,0)$, let $v \in V(H')$. Then, since $f_2(v)=f_2'(v) > 0$ and by Claim~\ref{claim_multiplicity}(a),  the pair $(H'',f'')=(H,v)/(f,2)$ is a hard pair of type (M) with  $f_1''(u) = d_{H''}(u) > 0$ for any $u \in V(H'')$. Since $f_1''(u)=f_1(u)=f_1(v)$ for any $u \in V(H'') \setminus \{z\}$ (as $f_1$ is constant in $V(H) \setminus \{z\}$), this implies that $f_1(v) > 0$ and $(H'',f^*)=(H,v)/(f,1)$ is a hard pair of type (M). However, it holds $f_1^*(z)=f_1(z) > 0$ and $f_2^*(u) = f_2(u) > 0$ for all $u \in V(H'') \setminus \{z\}$, which is impossible.

If $f(z) \neq (d_H(z),0,0,\ldots,0)$, there is an index $j \neq 1$ such that $f_j(z) > 0$. Since $E_H(z)$ contains only hyperedges, this implies that $(H,f)/(v,k)$ is not of type (M) for any $v \in V(H')$ and for any $k \in \{1,2,\ldots,p\}$ with $f_k(v) > 0$. Thus, after shrinking $H$ at any vertex, no hyperedge may remain. Nevertheless, since $n \geq 4$, this is impossible. This completes the proof.
\end{proof2}

\begin{claim}\label{claim_cycle}
The hard pair $(H',f')$ is not of type {\upshape (C)}.
\end{claim}

\begin{proof2}
Assume, to the contrary, that $(H',f')=(H,f)/(z,1)$ is of type (C) and, thus, $H=tC_{n-1}$ for some $t \geq 1, n \geq 6$ odd. Moreover, there are two indices $k \neq \ell$ from $\{1,2,\ldots,p\}$ such that $f'_k(v)=f'_\ell(v)=t$  and $f_j'(v)=0$ for all $j \in \{1,2,\ldots,p\} \setminus \{k,\ell\}$ and for all $v \in V(H')$. By symmetry, we may assume $k=2$ and $\ell \in \{1,3\}$. If $\ell=1$, then we obtain
\begin{itemize}
\item[$\circledast$] $f'(v)=(t,t,0,0,\ldots,0)$ and  $f(v)=(t + \mu_H(v,z), t, 0, 0, \ldots, 0)$
\end{itemize}
for all $v \in V(H')$. If $\ell =3$, it holds
\begin{itemize}
\item[$\circledcirc$]  $f'(v)=(0,t,t,0,0,\ldots,0)$ and $f(v)=(\mu_H(v,z),t,t,0,0,\ldots,0)$
\end{itemize}
for all $v \in V(H'$). Similar to the proof of Claim~\ref{claim_complete}, we distinguish between two cases.

\medskip
\case{A}{$E_H(z)$ contains an ordinary edge.} Then, the set $$N = \{v \in V(H') ~ | ~ vz \in E(H)\}$$ is non-empty. Since $f_2(v)=f_2'(v)=t > 0$ for all $v \in V(H')$ (by \eqref{eq_degreefirstcase}), this implies that $f_2(z) > 0$ (by Claim~\ref{claim_multiplicity}(b)) and so $(H',f'')=(H,f)/(z,2)$ is a hard pair of type (M) or (C). If $\circledast$ holds, then
$$f''(v)=(t + \mu_H(v,z), t - \mu_H(v,z), 0, 0, \ldots, 0)$$ for all $v \in V(H')$ and since $\mu_H(v,z) \geq 1$ for all $v \in N$ this implies that $(H',f'')$ is a bad pair of type (M). Then we conclude that $t - \mu_H(v,z)=0$ for all $v \in V(H')$ and so $\mu_H(v,z)=t$ for all $v \in V(H')$ and $N=V(H')$. If $\circledcirc$ holds, we have
$$f''(v)=(\mu_H(v,z),t,t,0,0,\ldots,0)$$ for all $v \in V(H')$ and again, since $\mu_H(v,z) \geq 1 $ for all $v \in N$, this implies that $(H',f'')$ is a bad pair of type (C). Hence, in both cases we have $\mu_H(v,z)=t$ for all $v \in V(H')$ and $N=V(H')$. Thus, we obtain
\begin{itemize}
\item[$\tilde{\circledast}$] $f(v)=(2t,t,0,0,\ldots,0)$ for all $v \in V(H')$ (if $\ell=1$), or
\item[$\tilde{\circledcirc}$] $f(v)=(t,t,t,0,0,\ldots,0)$ for all $v \in V(H')$ (if $\ell=3$).
\end{itemize}

Since $z$ is joined in $H$ to all other vertices by ordinary edges and since $H \div v \div z = H \div z \div v$ is a path (with multiple edges) and therefore connected for any $v \in V(H')$, $H \div v$ is a block for all $v \in V(H)$. As a consequence, for any vertex $v \in V(H')$, the hard pair $(H,f)/(v,1)$ must be of type (M). However, since for $H$ either $\tilde{\circledast}$ or $\tilde{\circledcirc}$ holds, it is easy to check that this is impossible.

\medskip
\case{B}{$E_H(z)$ contains only hyperedges.} As a consequence, $f(v)=f'(v)$ for all $v \in V(H')$.
We claim that $H$ admits an $f$-partition. To this end, let $e \in E(H)$ be an arbitrary hyperedge of $H$, and let $i_H(e)=\{x,y,z\}$  ($e$ must contain $z$ since $H \div z$ is a $tC_{n-1}$). Then, since $H'=H \div z$, the vertices $x$ and $y$ are adjacent in $H'$. Let $v_1,v_2,\ldots,v_{n-1}$ be a cyclic order of the vertices of $H'=tC_{n-1}$ with $x=v_1$ and $y=v_{n-1}$. Then, we define $H_2=H[ \{v_i \in V(H')~|~i \geq 1~\text{odd}\}]$, $H_\ell=H[\{v_i \in V(H')~|~i \geq 2~\text{even}\}],$ and $H_j=\ems$ for all $j \in \{1,2,\ldots,p\} \setminus \{2, \ell\}$. Since $d_{H_2}(v_1)=d_{H_2}(v_{n-1}) < t =f_2(v_1) = f_2(v_{n-1})$, $E(H_\ell)=\ems$, and $f_\ell(v)=t$ for all $v \in V(H')$ (see $\circledast$, $\circledcirc$), the  sequence $(H_1,H_2,\ldots,H_p)$ is an $f$-partition of $H-z$. As $H_1$ is an edgeless induced subhypergraph of $H \div z$, it follows that $d_{H_1+z}(z) = 0 < f_1(z)$, contradicting Proposition~\ref{prop_degree-condition}(c). This proves the claim.
\end{proof2}

\begin{claim}\label{claim_Hdiv z is block}
For every vertex $z\in V(H)$, $H\div z$ is not a block.
\end{claim}

\begin{proof2}
Suppose, to the contrary, that there exists a vertex $z$ such that $H \div z$ is a block. Let $j \in \{1,2,\ldots,p\}$ such that $f_j(z)>0$. Then, by Claim~\ref{claim_multiplicity}(a), $(H',f')=(H,f)/(z,j)$ is a hard pair and, since $H'=H\div z$ is a block, $(H',f')$ must be of type (M), (K), or (C). However, the three above claims imply that this is not possible.
\end{proof2}

\begin{claim}\label{claim_end-blocks}
For every vertex $z\in V(H)$, $H \div z$ has exactly two end-blocks.
\end{claim}

\begin{proof2}
Assume, to the contrary, that there is a vertex $z \in V(H)$ such that $H'=H \div z$ does not have exactly two end-blocks. By Claim~\ref{claim_Hdiv z is block} this implies that $H'$ has at least three end-blocks.

Let $T$ denote the block graph of $H'$, that is, the simple graph having vertex set $V(T)=\cB(H') \cup S$, where $S$ is the set of all separating vertices of $H'$, and edge set $E(T)=\{vB ~|~ v\in S, B\in \cB(H') \text{ and } v\in V(B)\}$. Note that $T$ is a tree with bipartition $(\cB(H'),S)$ and the end-blocks of $H'$ coincide with the leafs of $T$.
Since $H'$ has at least three end-blocks, $\Delta(T)\geq 3$. Let $B$ be an arbitrary end-block of $H'$. Since $B$ is a leaf of $T$ and $\Delta(T)\geq 3$, there is a unique vertex $x_B \in V(T)$ such that $x_B$ is the only vertex of degree at least $3$ in $T$ belonging to the subpath $P_B$ of $T$ between $x_B$ and $B$. Moreover, there exists a unique subtree $T_B$ of $T$ such that $T=T_B \cup P_B$ and $V(T_B) \cap V(P_B)= \{x_B\}$. Finally, there is a unique vertex $v_B\in S$ such that $v_B=x_B$ or $x_Bv_B$ is an edge of $P_B$.

Let $B_1,B_2,B_3$ be three distinct end-blocks of $H'$. For $i \in \{1,2,3\}$, let $v_i \in V(B_i)$ be the only separating vertex of $H'$ contained in $B_i$, let $u_i \in V(B_i) \setminus \{v_i\}$, let
$V_i$ be the set of all vertices contained in a block of $H'$ belonging to $T_{B_i}$, and let
$\tilde{B^i}=H[V_i \cup \{z\}]$. Since $H$ is a block, for each end-block $B$ of $H'$ there is an edge $e\in E_H(z)$ such that the vertex set $i_H(e)-\{z \}$ belongs to $B$ and contains a non-separating vertex of $H'$. Consequently, $\tilde{B^i}$ is a block contained in $H\div u_i$ as an induced subhypergraph and, therefore, there is a unique block $B^i$ of $H\div u_i$ containing $\tilde{B^i}$. Furthermore, let $(H^i,f^i)=(H,f)/(u_i,j_i)$ for some $j_i \in \{1,2,\ldots,p\}$ satisfying $f_{j_i}(u_i)>0$. Since $(H^i,f^i)$ is a hard pair (by Claim~\ref{claim_multiplicity}(a)), it follows from Proposition~\ref{prop_block-hardpair} that $(H^i,f^i)$ resulted from hard pairs $(B',f^i_{B'})$ (with $B' \in \mathcal{B}(H^i)$) of type (M), (K) and (C) by merging them appropriately ($H^i$ has at least two blocks by Claim~\ref{claim_Hdiv z is block}). Note that for $i \neq j$ from $\{1,2,3\}$, the vertices $u_i$ and $u_j$ are not adjacent in $H \div z$ and, therefore, not adjacent in $H$. As a consequence, the hard pair $(B^i,f^i_{B^i})$ cannot be of type (K) for $i \in \{1,2,3\}$. Furthermore, for $i \in \{1,2,3\}$ and for $j \in \{1,2,3\}\setminus\{i\}$ the vertex $u_i$ is not a separating vertex of $H^j$ contained in $B^j$ and, together with Proposition~\ref{prop_block-hardpair}, we conclude $$\circledast \quad f^j_{B^j}(u_i)=f^j(u_i)=f(u_i).$$ In the following, we regard the hard pairs $(B^1,f^1_{B^1}), (B^2,f^2_{B^2})$ and $(B^3,f^3_{B^3})$.

\medskip
\case{A}{One of the three hard pairs, say $(B^1,f^1_{B^1})$ is of type {\upshape (C)}.} Then, $B^1=tC_m$ for some $t \geq 1, m \geq 5$ odd and by symmetry we may assume $f^1_{B^1}(v)=(t,t,0,\ldots,0)$ for all $v \in V(B)$. Since $\tilde{B^1}$ contains no separating vertex and as $B^1$ is a $tC_m$, $\tilde{B^1}=B^1$. Furthermore, this implies that $z$ is joined to $u_2$ and $u_3$ by ordinary edges in $H$ (since $u_1$ and $u_i$ are not adjacent in $H \div z$ for $i \in \{2,3\}$), that $V(B_i)=\{u_i,v_i\}$ for $i \in \{2,3\}$, and that $P^1=T_{B_1}$ is a path and each block on $P^1$ is a $tK_2$. Since $f^3_{B^3}(u_2)=f(u_2)=f^1_{B^1}(u_2)=(t,t,0,0,\ldots,0)$ (by $\circledast$) and since $(B^3,f^3_{B^3})$ is a hard pair (not of type (K)) with $u_2 \in V(B^3)$ not being a separating vertex of $H^3$, we obtain that $(B^3,f^3_{B^3})$ must be a hard pair of type (C), too, and that $P^3=T_{B_3}$ is a path and each block on $P^3$ is a $tK_2$. Analoguesly we can show that $(B^2,f^2_{B^2})$ is a hard pair of type (C) and that each block of the path  $P^2=T_{B_2}$ is a $tK_2$. Since $u_1,u_2$ and $u_3$ are not pairwise adjacent in $H$, this implies that $B^1=tC_m$ contains exactly one separating vertex $v_B$ of $H^1$ and that $H$ is the union of the three (multi-)cycles $B^1,B^2,B^3$ with $V(B^1) \cap V(B^2) \cap V(B^3) = \{z,v_B\}$ and $u_i \not \in V(B^i)$. Let $\ell_i$ be the length of the $(z,v_B)$-(multi-)path  in $H$ containing the ordinary edge $zu_i$. Then,
\begin{align*}
|B^1| & =  \ell_2 + \ell_3, \\
|B^2| & =  \ell_1 + \ell_3, \text{ and}\\
|B^3| & =  \ell_1 + \ell_2.
\end{align*}
However, since $|B^i|$ is odd for $i \in \{1,2,3\}$, we obtain $$\ell_1 + \ell _2 \equiv \ell_1 + \ell_3 \equiv \ell_2 + \ell_3 \equiv 1 \text{ (mod}~p\text{)},$$ which is impossible.

\medskip
\case{B}{All three hard-pairs $(B^i,f^i_{B^i})$ $(i \in \{1,2,3\})$ are of type {\upshape(M)}.}
By symmetry, we may assume that $f^1_{B^1}(v)=(d_{B^1}(v),0,0,\ldots,0)$ for all $v \in V(B^1)$. Since the vertex $u_2$ is contained in $B^1$ and $B^3$, and since $f^1_{B^1}(u_2)=f^1(u_2)=f(u_2)=f^3_{B^3}(u_2)$ (by $\circledast$), we conclude that $f^3_{B^3}(v)=(d_{B^3}(v),0,0,\ldots,0)$ for all $v \in V(B^3)$. Analoguesly, regarding $u_3 \in V(B^2)$, we conclude that $f^2_{B^2}(v)=(d_{B^2}(v),0,0,\ldots,0)$ for all $v \in V(B^2$).
As $(H^i,f^i)=(H,f)/(u_i,j_i)$, this implies that $j_i=1$ for $i \in \{1,2,3\}$. Furthermore, since $z$ clearly is a non-separating vertex of $H^i$ contained in $B^i$ for all $i \in \{1,2,3\}$ and as $j_i=1$ for $i \in \{1,2,3\}$, it follows $f(z)=(d_H(z),0,0,\ldots,0)$. Let $v \in V(H) \setminus \{z\}$ be an arbitrary vertex. We claim that $f(v)=(d_H(v),0,0,\ldots,0)$. Assume this is false. Then, by symmetry, $f_2(v)>0$. Clearly, $v$ belongs to $\tilde{B^i}$
for some $i \in \{1,2,3\}$ , say $v$ belongs to $\tilde{B^1}$. Note that $\tilde{B^1}$ is an induced subhypergraph of $B^1\in \cB(H^i)$. If $v$ is a non-separating vertex of $H^1$, then it holds $f(v)=(d_H(v),0,0,\ldots,0)$ (as $j_1=1$, as $f^1(v)=f^1_{B^1}(v)=(d_{B^1}(v),0,0,\ldots,0)$ and by \eqref{eq_prop2.5}), contradicting our assumption. Otherwise, $v=v_{B_1}$ is a separating vertex of $H^1$ and so $\tilde{B^1}=B^1$. As $(H^1,f^1)$ results from merging hard pairs (by Proposition~\ref{prop_block-hardpair}) and as $f_2(v)>0$, this implies that $v$ is contained in a block $B'$ of $H^1$ with $V(B') \cap V(B^1)=\{v\}$ and $f_2(w)>0$ for all $w \in V(B')$. Let $v' \in V(B') \setminus \{v\}$. Then, $v'$ is contained in a block belonging to the subpath of the block graph $T$ between $v_{B_1}$ and $B_1$. But then, $v'$ is a non-separating vertex of $H^j$ contained in $\tilde{B^j}$ for $j\in \{2,3\}$. This however implies that $f_{B^2}(v')=(d_{B^2}(v'),0,0,\ldots,0)$, a contradiction. Hence, the claim is proved and, thus, $(H,f)$ is a hard pair of type (M), contradicting (2).
\end{proof2}

\begin{claim} \label{claim_blocks}
There exists a sequence $B_1,B_2,\ldots,B_\ell$ of $\ell \geq 4$ induced subhypergraphs of $H$ and a sequence $u_0, u_1, \ldots, u_{\ell - 1}$  of distinct vertices of $H$  such that the following statements hold.
\begin{itemize}
\item[\upshape(a)] $B_i=t_iK_2$ for $i \in \{2,3,\ldots,\ell\}$ and some $t_i \geq 1$, $B_1$ has no separating vertex, and $|B_1| \geq 2$.
\item[\upshape(b)] $H=B_1 \cup B_2 \cup \ldots \cup B_\ell$, $V(B_i) \cap V(B_{i+1})=\{u_{i}\}$ for $i \in \{1,2,\ldots, \ell-1\},$ and $V(B_1) \cap V(B_\ell) = \{u_0\}$.
\end{itemize}
\end{claim}

\begin{figure}[thbp]
\centering \resizebox{.8\linewidth}{!}{
\begin{tikzpicture} [node distance=1cm, bend angle=20,
vertex/.style={circle,minimum size=3mm,very thick, draw=black, fill=black, inner sep=0mm}, information text/.style={fill=red!10,inner sep=1ex, font=\Large}, help lines/.style={-,color=black}]

\draw (-3,1) ellipse [x radius=2cm, y radius=1cm, rotate=-20];

\node[vertex, label={[font=\Large]90: $u_{\ell - 1}$}] (v) at (3,5) {};
\node[vertex,  label={[font=\Large]230: $u_1$}] (u1) at (-1.5,.5) {};
\node[vertex, label={[font=\Large]270: $u_2$}] (u2) at (1.5, 0) {};
\node[vertex, label={[font=\Large]270: $u_{\ell-3}$}] (u3) at (4.5, 0) {};
\node[vertex, label={[font=\Large]270: $u_{\ell-2}$}] (u4) at (7.5,.5) {};
\node [font=\Large] at (2.9,0) {$\ldots$};
\node[vertex, label={[font=\Large]270: $u_0$}] (u0) at (-4,1.5) {};

\node [font=\Large] at (-3.5,-.5) {$B_1$};
\node [font=\Large] at (-.3, -.5) {$B_2$};
\node [font=\Large] at (6.3, -.5) {$B_{\ell-2}$};
\node [font=\Large] at (7,3.25) {$B_{\ell-1}$};
\node [font=\Large] at (-1,4.5) {$B_{\ell}$};

%
%

\path[-]
(v) edge [help lines, bend angle=20, bend right] node [sloped, below] {$t_\ell$} (u0)
		edge [help lines, bend angle=20, bend left] node [sloped, below] {$t_{\ell-1}$}(u4)
(u1) edge[help lines, bend angle=5, bend right] node [sloped, above] {$t_{2}$}
  (u2)

(u3) edge[help lines, bend angle=5, bend right] node [sloped, above] {$t_{\ell-2}$}
  (u4);

\end{tikzpicture}
}
\caption{The structure of $H$.}
\end{figure}
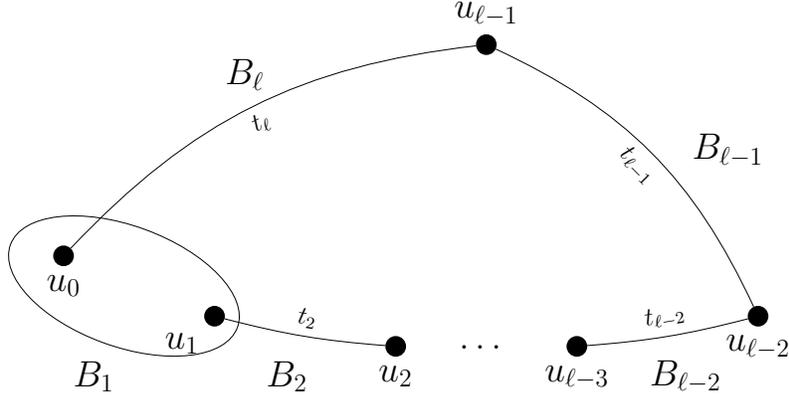

\begin{proof2}
Let $z$ be an arbitrary vertex of $H$. Then, $H \div z$ has exactly two end-blocks (by Claim~\ref{claim_end-blocks}) and, therefore, there is a uniquely determined sequence $B_1,B_2,\ldots,B_k$ of $k \geq 2$ blocks of $H \div z$ and a sequence $u_1,u_2,\ldots,u_{k-1}$ of distinct vertices such that $V(B_i) \cap V(B_{i+1})=\{u_i\}$ for all $i \in \{1,2,\ldots,k-1\}$ and $H \div z = B_1 \cup B_2 \cup \ldots \cup B_k$. In particular, $B_1$ and $B_k$ are the end-blocks of $H \div z$. Let $b_z=\max \{|B_i|~|~ i \in \{1,2,\ldots,k\}\}$. Among all vertices $z$ of $H$ we may choose one for which $b_z$ is maximum. Let $B_j$ be a block of $H \div z$ with $|B_j|=b_z$. Since $H$ is a block, there are vertices $u_0 \in V(B_1)$  and $u_k \in V(B_k)$ which are non-separating vertices of  $H \div z$ and adjacent to $z$. Assume that there is an index $i \neq j$ from the set $\{1,2,\ldots,k\}$ such that $B_i$ contains a non-separating vertex $v$ of $H \div z$ different from $u_0$ and $u_k$. Then, it follows that $H \div v$ has a block $B$ containing $B_j$ as well as $z$ and, thus, $|B| > |B_j|=b_z$, a contradiction. As a consequence,  for each index $i$ from the non-empty set $ \{1,2,\ldots,k\} \setminus \{j\}$, there exists an integer $t_i \geq 1$ such that $B_i=t_iK_2$. To complete the proof, all we need to show is that $z$ is not adjacent to any vertex besides $u_0$ and $u_k$. By symmetry, we may assume that $j \neq 1$ and, thus, $B_1$ is a $t_1K_2$ for some $t_1 \geq 1$ and $V(B_1)=\{u_0,u_1\}$. If there was a hyperedge $e$ with $i_H(e)=\{z,u_0,u_1\}$, then clearly $H \div u_0$ would still be a block, which is impossible. Thus, $u_0$ is adjacent only to $z$ and $u_1$ and not contained in any hyperedge. As a consequence, if $z$ is adjacent to any vertex from $V(B_1) \cup V(B_2) \cup \ldots \cup V(B_j)$ except from $u_0$ and $u_k$, then $H \div u_0$ has a block $B$ that contains $B_j$ as well as $z$, giving a contradiction to the maximality of $|B_j|$. Similar, if $z$ is adjacent to any vertex from $V(B_{j+1}) \cup V(B_{j+1}) \cup \ldots \cup V(B_k)$ except $u_k$ (implying $j \neq k$), by a similar argumentation we conclude that $H \div u_k$ has a block $B$ containing $B_j$ as well as $z$, again contradicting the maximality of $|B_j|$. Consequently, $z$ is only adjacent to $u_0$ and $u_k$. By setting $B_{k+1}=H[\{z,u_k\}$, $B_{k+2}=H[\{z,u_0\}]$, $u_{k+1}=z$ and by shifting the block-sequence and the vertex-sequence we obtain the required statement.
\end{proof2}

To conclude the proof, we show that $(H,f)$ is a hard pair, giving a contradiction to statement (2). Consider the sequences $B_1,B_2,\ldots,B_\ell$ and $u_0,u_1,\ldots,u_{\ell-1}$ as described in Claim~\ref{claim_blocks}. For technical reasons, let $u_\ell=u_0$.  Let $j \in \{1,2,\ldots,p\}$ such that $f_j(u_{\ell - 1})>0$, and consider the hard pair $(H',f')=(H,f)/(u_{\ell-1},j)$. Then, $B_1$ is a block of $H'=H \div u_{\ell - 1}$. Let $f'_{B_1}$ be defined as in Proposition~\ref{prop_block-hardpair}. We claim that $(B_1,f'_{B_1})$ is of type (M). If $(B_1, f'_{B_1})$ is of type (K), then $|B_1| \geq 3$ and there is a vertex $v \in V(B_1)$ such that $H \div v$ is a block, which is impossible. If $(B_1,f'_{B_1})$ is of type (C), say $B_1=tC_m$ with $m \geq 5$ odd, by symmetry we may assume that $f'_{B_1}(v)=(t,t,0,0,\ldots,0)$ for all $v \in V(B_1)$. Let $P_1$ and $P_2$ be the two disjoint $(u_0,u_1)$-(multi-)paths in $B_1$ and let $\ell_i$ be the length of $P_i$ for $i \in \{1,2\}$.
Furthermore, for $i \in \{1,2\}$, let $v_i$ be a vertex of $V(P_i) \setminus \{u_0,u_1\}$. Then, by regarding the hard pair $(H^1,f^1)=(H,f)/(v_1,1)$, we conclude that $B^2=H[V(P_2) \cup V(B_2) \cup V(B_3) \cup \ldots \cup V(B_{\ell})]$ is a block of $H^1$. Since $f^1(v_2)=f(v_2)=f'_{B_1}(v_2)=(t,t,0,0,\ldots,0)$, it follows that $(B^2, f^1_{B^2})$ is of type (C) and so $B^2=tC_{m'}$ for some $m' \geq 5$ odd. By a similar argumentation we obtain that, regarding $(H^2,f^2)=(H,f)/(v_2,1)$, $B^1=H[V(P_1) \cup V(B_2) \cup V(B_3) \cup \ldots \cup V(B_\ell)]$ must be a $tC_{m''}$ of odd length, as well. However, this implies that $\ell_1 + \ell_2 \equiv 0~\text{(mod 2)}$ and so $B_1=tC_m$ with $m$ even, which is impossible. As a consequence, $(B_1,f_{B_1}')$ is of type (M) as claimed and so there is an index $j' \in \{1,2,\ldots,p\}$ such that $f_{j'}(v)=d_{B_1}(v)=d_H(v)$ and $f_k(v)=0$ for $k \in \{1,2,\ldots,p\} \setminus \{j'\}$ and for all $v \in V(B_1) \setminus \{u_0,u_1\}$. Moreover, it holds $f_{j'}(v) \geq d_{B_1}(v)$ for all $v \in \{u_0,u_1\}$.

Recall that for $i \in \{2,3,\ldots, \ell \}$, we have $B_i=t_iK_2$ and $V(B_i)=\{u_{i-1},u_i\}$.
Furthermore, recall that there is an index $j \in \{1,2,\ldots,p\}$ such that $f_j(u_{\ell - 1}) > 0$. By symmetry, we may assume that $j=1$. We claim that either
\begin{itemize}
\item[$\circledast$]$f(u_i)=(t_{i}+t_{i+1},0,0,\ldots,0)$ for all $i \in \{2,3,\ldots,\ell-1\}$ , or
\item[$\circledcirc$] $f(u_i)=(t,t,0,0,\ldots,0)$ for all $i \in \{2,3,\ldots, \ell - 1\}$ (except for symmetry) and $B_i=tK_2$
for all $i \in \{2,3,\ldots,\ell\}$.
\end{itemize}
Since $f_1(u_{\ell - 1}) > 0$, by repeated application of Claim~\ref{claim_multiplicity}(b) we conclude $f_1(u_i) \geq \max\{t_{i},t_{i+1}\}$ for $i \in \{2,3,\ldots,\ell-1\}$. If there exists an index $k \neq 1$, say $k=2$ (by symmetry) such that $f_2(u_i) > 0$ for some $i \in \{2,3,\ldots,\ell-1\}$, then, similarly to above, we get $f_2(u_i) \geq \max\{t_{i},t_{i+1}\}$ for $i \in \{2,3,\ldots,\ell-1\}$. By \eqref{eq_prop2.5}, this implies $t_{i}=t_{i+1}=t$ as well as  $f(u_i)=(t,t,0,0,\ldots,0)$ for some $t \geq 1$ and for all $i \in \{2,3,\ldots, \ell - 1\}$, and so $\circledcirc$ holds. If $f_k(u_i)=0$ for all $k \in \{2,3,\ldots,p\}$ and all $i \in \{2,3,\ldots,\ell-1\}$, equation \eqref{eq_prop2.5} implies that $f(u_i)=(t_{i}+t_{i+1},0,0,\ldots,0)$ for all $i \in \{2,3,\ldots,\ell-1\}$ and $\circledast$ holds.

If $\circledast$ is satisfied, then by Claim~\ref{claim_multiplicity}(b) it holds $f(v)=(d_H(v),0,0,\ldots,0)$ for $v \in \{u_0,u_1\}$ and hence $j'=1$ and $(H,f)$ is a hard pair of type (M), contradicting (2).

Thus, it remains to consider the case that $\circledcirc$ holds. If $|B_1|=2$, then $B_1=t_1K_2$ for some $t_1 \geq 1$. Then, again we conclude $f(u_0)=f(u_1)=(t,t,0,0,\ldots,0)$ and so $H=tC_n$ and $f(v)=(t,t,0,0,\ldots,0)$ for all $v \in V(H)$. Furthermore, $n$ must be odd since otherwise $H$ would clearly admit an $f$-partition. Consequently, $(H,f)$ is of type (C), which contradicts (2).

Finally, assume that $|B_1| \geq 3$. Then, there is a vertex $z \in V(B_1)$ different from $u_0$ and $u_1$ and $f_{j'}(z)=f'_{j'}(z)=d_{B_1}(z)=d_H(z)$ and
$f_k(z)=0$ for $k \in \{1,2,\ldots,p\} \setminus \{j'\}$. As $f(u_{\ell-1})=(t,t,0,0,\ldots,0)$, it follows from Claim~\ref{claim_multiplicity}(b) that $f_1(u_0)\geq t > 0$ and $f_2(u_0)\geq t >0$. Since $(H',f')=(H,f)/(u_{\ell-1},1)$ and since $(B_1, f_{B_1}')$ is a hard pair of type (M), it must hold $f_1(u_0)=t$ and $j'=2$. Therefore, we have $f(z)=f'(z)=(0,d_H(z),0,0,\ldots,0)$. Moreover, as $f_2(u_{\ell - 1})=t > 0$, $(H',f'')=(H,f)/(u_{\ell-1},2)$ is a hard pair, too, and $(B_1,f_{B_1}'')$ is a hard pair of type (M). Consequently, it must hold $f_2(u_0)=t$ and $f_1''(v) > 0$ for all $v \in V(B_1)$. However, as $f_1''(z)=f_1(z)=f_1'(z)=0$, this is impossible. This completes the proof.
\end{proof}

\section{Applications of Theorem \ref{Theorem:Hauptsatz}}
In this section, some applications of Theorem~\ref{Theorem:Hauptsatz} are presented. Those are guided by the depiction of Borodin, Kostochka and Toft in \cite{BorKosToft}.

\subsection{Brooks' Theorem for list-colorings of hypergraphs}
Recall from the introduction that the chromatic number, respectively list-chromatic number of a hypergraph $H$ is always less or equal to the coloring number of $H$. In particular, it holds
$$\chi(H) \leq \chi^\ell(H) \leq \text{col}(H) \leq \Delta(H) + 1.$$
This inequality naturally raises the question, in which cases $\chi(H)=\Delta(H) + 1$ holds. For simple graphs, the answer was given by Brooks \cite{Brooks} in 1941. His famous theorem states that complete graphs and odd cycles are the only connected graphs, for which the chromatic number is equal to the maximum degree plus one. For list-colorings, the solution was found by Erd\H os, Rubin and Taylor \cite{erd} and, independently, by Vizing \cite{Vizing}. They proved the following.

\begin{theorem}{\bf(Erd\H os, Rubin and Taylor)}
Let $G$ be a connected simple graph and let $L$ be a list-assignment satisfying $|L(v)| \geq d_G(v)$ for all $v \in V(G)$. If $G$ does not admit a proper $L$-coloring, then $|L(v)|=d_G(v)$ for all $v \in V(H)$ and each block of $G$ is either a complete graph or an odd cycle. As a consequence, $\chi^\ell(G) \leq \Delta(G) + 1$ and equality holds if and only if $G$ is a complete graph or an odd cycle.
\end{theorem}

It shows that those two theorems can be extended to hypergraphs, as well. An analogue to Brooks' Theorem was given by Jones \cite{Jones} in 1975. Brooks' Theorem for list-colorings of hypergraphs was obtained by Kostochka, Stiebitz and Wirth \cite{KoStiWi} in 1995.

\begin{theorem}{\bf(Kostochka, Stiebitz and Wirth)} \label{Theorem:KoStiWi}
Let $H$ be a connected simple hypergraph and let $L$ be a list-assignment satisfying $|L(v)| \geq d_H(v)$ for all $v \in V(H)$. If $H$ does not admit a proper $L$-coloring, then it holds $|L(v)|=d_H(v)$ for all $v \in V(H)$ and each block $B$ of $H$ is either a complete graph, an odd cycle, or $B$ has just one edge. As a consequence, $\chi^\ell(H) \leq \Delta(H) + 1$ and equality holds if and only if $H$ is either a complete graph, an odd cycle, or if $H$ contains just one edge.
\end{theorem}

How can we conclude the above theorem from Theorem~\ref{Theorem:Hauptsatz}? To this end, let $H$ be a simple hypergraph and let $L$ be a list-assignment of $H$ with a set $C$ of $p$ colors. By renaming the colors from $C$ we may assume $C=\{1,2,\ldots,p\}$. Let $f \in \mathcal{V}_p(H)$ such that $f_i(v)=1$ if $i \in L(v)$ and $f_i(v)=0$ otherwise for all $i \in \{1,2,\ldots,p\}$ and for all $v \in V(H)$. Then, for any $L$-coloring $\varphi$ of $G$, the sequence $(H_1,H_2,\ldots,H_p)$ with $H_i = H[\varphi^{-1}(i)]$ for all $i \in \{1,2,\ldots,p\}$ is an $f$-partition of $H$. Conversely, if $(H_1,H_2,\ldots,H_p)$ is an $f$-partition of $H$, then setting $\varphi(v)=i$ if $v \in V(H_i)$ leads to an $L$-coloring of $H$, since $H_i$ is edgeless by construction. As a consequence, $H$ admits an $L$-coloring if and only if $H$ is $f$-partitionable.

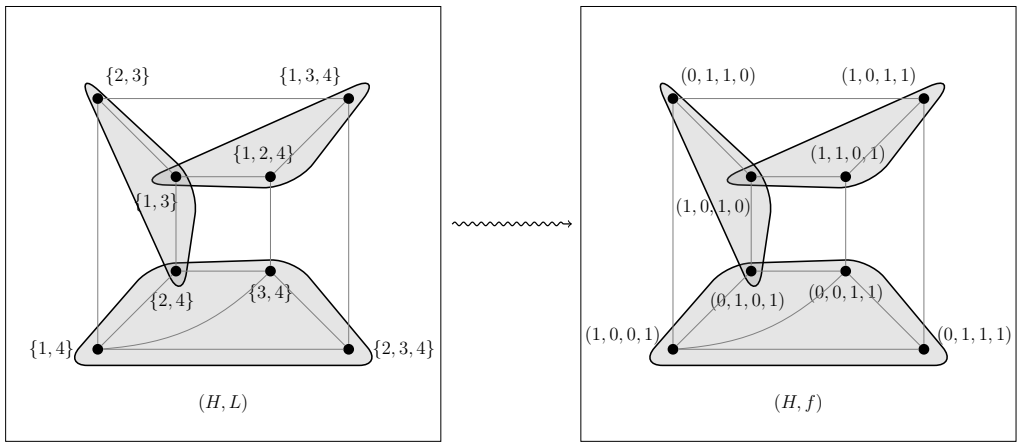
\begin{figure}[thbp]
\centering \resizebox{\linewidth}{!}{
%
\begin{tikzpicture}
[node distance=1cm, bend angle=20,
vertex/.style={circle,minimum size=2mm,very thick, draw=black, fill=black, inner sep=0mm}, information text/.style={fill=red!10,inner sep=1ex, font=\Large}]

\begin{scope}[xshift=-5cm]

\node[draw=none,minimum size=3cm,regular polygon,regular polygon sides=4] (a) {};
\foreach \x in {1,2,3,4}
\node[vertex] (v\x) at (a.corner \x) {};

\node[draw=none,minimum size=8cm,regular polygon,regular polygon sides=4] (b) {};
\foreach \x in {1,2,3,4}
\node[vertex] (u\x) at (b.corner \x) {};

\node at (v1) [label={[xshift=20pt]93:$\{1,2,4\}$}] {};
\node at (v2) [label={[xshift=10pt, yshift=-5pt]200:$\{1,3\}$}] {};
\node at (v3) [label={[xshift=20pt, yshift=-8pt]200:$\{2,4\}$}] {};
\node at (v4) [label={270:$\{3,4\}$}] {};

\node at (u1) [label={100:$\{1,3,4\}$}] {};
\node at (u2) [label={80:$\{2,3\}$}] {};
\node at (u3) [label={[xshift=-7pt]180:$\{1,4\}$}] {};
\node at (u4) [label={[xshift=7pt]0:$\{2,3,4\}$}] {};

\begin{pgfonlayer}{background}
\filldraw[fill=black!50!, fill opacity=.2,rounded corners=20pt, line width=1pt, draw=black] (-2,0.9)--(1.6,0.8)--(3.6,3.4)--cycle;
		
\filldraw[fill=black!50!, fill opacity=.2,rounded corners=15pt, line width=1pt, draw=black] (-1.6,-0.9)--(1.6,-0.8)--(3.6,-3.2)--(-3.6,-3.2)--cycle;

\filldraw[fill=black!50!, fill opacity=.2,rounded corners=20pt, line width=1pt, draw=black, xshift=-1pt] (-0.5,0.9)--(-3.3,3.5)--(-0.9,-1.8)--cycle;		
\end{pgfonlayer}		

\path[-]
		(v1)
		edge[help lines] (v2)
		edge[help lines] (u1)
		(v2)
		edge[help lines] (u2)
		edge[help lines] (v3)
		(v3)
		edge[help lines] (u3)
		edge[help lines] (v4)
		(v4)
		edge[help lines, bend left] (u3)
		edge[help lines] (u4)
		edge[help lines] (v1)
		(u1)
		edge[help lines] (u2)
		(u2)
		edge[help lines] (u3)
		(u3)
		edge[help lines] node [name=l1, midway, sloped, below=2pt] {} (u4)
		(u4)
		edge[help lines] (u1);
%
%

\node at (l1) [yshift=-1cm] {$(H,L)$};

\node (h1) at (u1) [xshift=1.8cm, yshift=1.8cm] {};
\node (h2) at (u2) [xshift=-1.8cm, yshift=1.8cm] {};
\node (h3) at (u3) [xshift=-1.8cm, yshift=-1.8cm] {};
\node (h4) at (u1) [xshift=1.8cm, yshift=-1.8cm] {};

\node (l1) [draw=black, rectangle, fit=(h1)(h2)(h3)(h4)] {};
\node (r1) at (l1.east) {};

\end{scope}

\begin{scope}[xshift=8cm]
\node[draw=none,minimum size=3cm,regular polygon,regular polygon sides=4] (a) {};
\foreach \x in {1,2,3,4}
\node[vertex] (v\x) at (a.corner \x) {};

\node[draw=none,minimum size=8cm,regular polygon,regular polygon sides=4] (b) {};
\foreach \x in {1,2,3,4}
\node[vertex] (u\x) at (b.corner \x) {};

\node at (v1) [label={[xshift=30pt]93:$(1,1,0,1)$}] {};
\node at (v2) [label={[xshift=8pt, yshift=-8pt]200:$(1,0,1,0)$}] {};
\node at (v3) [label={[xshift=30pt, yshift=-8pt]200:$(0,1,0,1)$}] {};
\node at (v4) [label={270:$(0,0,1,1)$}] {};

\node at (u1) [label={100:$(1,0,1,1)$}] {};
\node at (u2) [label={80:$(0,1,1,0)$}] {};
\node at (u3) [label={[xshift=0pt, yshift=9pt]180:$(1,0,0,1)$}] {};
\node at (u4) [label={[xshift=0pt, yshift=9pt]0:$(0,1,1,1)$}] {};

\begin{pgfonlayer}{background}
\filldraw[fill=black!50!, fill opacity=.2,rounded corners=20pt, line width=1pt, draw=black] (-2,0.9)--(1.6,0.8)--(3.6,3.4)--cycle;
		
\filldraw[fill=black!50!, fill opacity=.2,rounded corners=15pt, line width=1pt, draw=black] (-1.6,-0.9)--(1.6,-0.8)--(3.6,-3.2)--(-3.6,-3.2)--cycle;

\filldraw[fill=black!50!, fill opacity=.2,rounded corners=20pt, line width=1pt, draw=black, xshift=-1pt] (-0.5,0.9)--(-3.3,3.5)--(-0.9,-1.8)--cycle;		
\end{pgfonlayer}

\path[-]
		(v1)
		edge[help lines] (v2)
		edge[help lines] (u1)
		(v2)
		edge[help lines] (u2)
		edge[help lines] (v3)
		(v3)
		edge[help lines] (u3)
		edge[help lines] (v4)
		(v4)
		edge[help lines, bend left] (u3)
		edge[help lines] (u4)
		edge[help lines] (v1)
		(u1)
		edge[help lines] (u2)
		(u2)
		edge[help lines] (u3)
		(u3)
		edge[help lines] node [name=l1, midway, sloped, below=2pt] {} (u4)
		(u4)
		edge[help lines] (u1);

\node at (l1) [yshift=-1cm] {$(H,f)$};

\node (h1) at (u1) [xshift=1.8cm, yshift=1.8cm] {};
\node (h2) at (u2) [xshift=-1.8cm, yshift=1.8cm] {};
\node (h3) at (u3) [xshift=-1.8cm, yshift=-1.8cm] {};
\node (h4) at (u1) [xshift=1.8cm, yshift=-1.8cm] {};

\node (l2) [draw=black, rectangle, fit=(h1)(h2)(h3)(h4)] {};
\node (r2) at (l2.west) {};

\end{scope}

\draw [shorten >=.5mm,-to,thick,decorate,
decoration={snake,amplitude=.4mm,segment length=2mm,
pre=moveto,pre length=1mm,post length=2mm}]
(r1) -- (r2); 

%
%

%
%
%
%
%
\end{tikzpicture}
}
\caption{Transforming a list-assignment $L$ into a function $f$.}
\end{figure}

Thus, as $|L(v)|\geq d_H(v)$ and so $f_1(v) + f_2(v) + \ldots + f_p(v) \geq d_H(v)$ for all $v \in V(H)$, if $H$ does not admit an $L$-coloring, Theorem~\ref{Theorem:Hauptsatz} implies that $(H,f)$ is a hard pair; therefore, each block from $H$ is of type (M), (K), or (C), and we can easily conclude the first part of Theorem~\ref{Theorem:KoStiWi}. In order to deduce the second part of the theorem, we argue as follows. If $\chi^\ell(H) = \Delta(H) + 1$, then there is a list-assignment $L$ satisfying $|L(v)|=\Delta(H)$ for all $v \in V(H)$ such that $H$ does not admit a proper $L$-coloring. Consequently, it must hold that $d_H(v)= |L(v)|=\Delta(H)$ for all $v \in V(H)$, and so $H$ is $\Delta(H)$-regular. Moreover, the first part of Theorem~\ref{Theorem:KoStiWi} also implies that each block from $H$ is a complete graph, an odd cycle or contains just one edge. But then, as $H$ is $\Delta(H)$-regular, $H$ can only consist of exactly one block, and we are done. On the other hand, if $H$ is a complete graph, an odd cycle, or if $|E(H)|=1$, then it easy to see that $\chi^\ell(H)=\Delta + 1$.

\subsection{Additional degree constraints}
Borodin \cite{Bor} and, independently, Bollob\'as and Manvel \cite{BolMan}, proved another extension of Brooks' Theorem for the class of simple ordinary graphs.

\begin{theorem}{\bf (Borodin/Bollob\'as and Manvel)} \label{Theorem:BorBolMan}
Let $G$ be a connected simple graph with maximum degree $\Delta  \geq 3$ different from $K_{\Delta + 1}$. Let also $k_1,k_2,\ldots,k_p$ be positive integers, $p \geq 2$, such that $k_1 + k_2 + \ldots + k_p \geq \Delta$. Then, there is a $p$-partition $(G_1,G_2,\ldots,G_p)$ of $G$ such that $\text{\upshape col}(G) \leq k_i$ whenever $1 \leq i \leq p$.
\end{theorem}

Clearly, by setting $k_1 = k_2 = \ldots = k_p = 1$ one can immediately deduce Brooks' Theorem from the above theorem. However, Borodin \cite{Bor2} generalized Theorem~\ref{Theorem:BorBolMan} even further with the help of a simple argument. Bollob\'as and Manvel \cite{BolMan} proved  the same extension independently.

\begin{theorem}{\bf (Borodin/Bollob\'as and Manvel)} \label{Theorem:col_partition}
Let $G$ be a connected simple graph with maximum degree $\Delta \geq 3$ different from $K_{\Delta + 1}$. Let also $k_1,k_2,\ldots,k_p$ be positive integers, $p \geq 2$, such that $$k_1 + k_2 + \ldots + k_p \geq \Delta.$$ Then, there is a $p$-partition $(G_1,G_2,\ldots,G_p)$ of $G$ satisfying $\text{\upshape col}(G) \leq k_i$ and $\Delta(G_i) \leq k_i$ whenever $1 \leq i \leq p$.
\end{theorem}

It shows that it is possible to prove a similar result for arbitrary hypergraphs.

\begin{stheorem}[\ref{Theorem:col_partition}]
\label{Theorem:col_partition_hypercase}
Let $H$ be a connected hypergraph having maximum degree $\Delta(H)=\Delta \geq 1$ that is not a $tK_n$ for some $t,n \geq 1$ and not a $tC_n$ for $t \geq 1, n \geq 3$ odd. Let also $k_1,k_2, \ldots, k_p$ be positive integers, $p \geq 2$, such that $k_1 + k_2 + \ldots + k_p \geq \Delta$. Then, there is a $p$-partition $(H_1,H_2,\ldots,H_p)$ of $H$ such that $\text{\upshape col}(H_i) \leq k_i$ and $\Delta(H_i) \leq k_i$ whenever $1 \leq i \leq p$.
\end{stheorem}

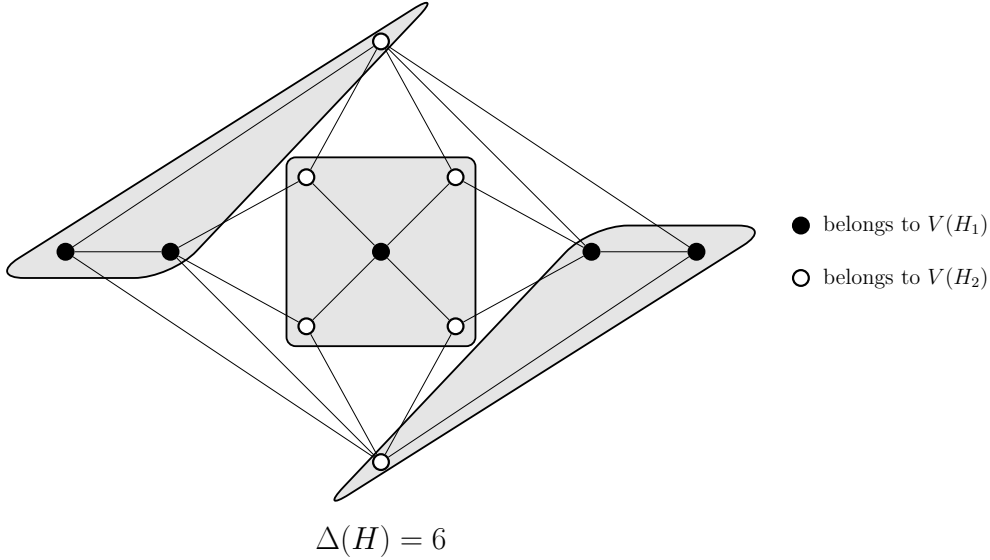
\begin{figure}[thbp]
\centering \resizebox{\linewidth}{!}{

\begin{tikzpicture} [node distance=1cm, bend angle=20,
vertex/.style={circle,minimum size=3mm,very thick, draw=black, fill=white, inner sep=0mm}, information text/.style={fill=red!10,inner sep=1ex, font=\Large}, help lines/.style={-,color=black}]

\node[vertex, fill=black] (v1) {};

\node[draw=none,minimum size=4cm,regular polygon,regular polygon sides=4] (b) {};
\foreach \x in {1,2,3,4}
\node[vertex] (w\x) at (b.corner \x) {};

\node[draw=none,minimum size=8cm,regular polygon,regular polygon sides=4, rotate=45] (c) {};
\foreach \x in {2,4}
\node[vertex, fill=black] (u\x) at (c.corner \x) {};
\node[vertex] (u3) at (c.corner 3){};
\node[vertex] (u1) at (c.corner 1){};

\node[vertex, fill=black] (v2) at (u2) [xshift=-2cm] {};
\node[vertex, fill=black] (v3) at (u4) [xshift=2cm] {};

\begin{scope}[xshift=8cm, yshift=.5cm]
\node[vertex, fill=black] (r1) {};
\node[vertex, fill=white, yshift=-1cm](r2){};
\node (t1) at (r1) [xshift=2cm] {belongs to $V(H_1)$};
\node (t2) at (r2) [xshift=2cm] {belongs to $V(H_2)$};
\end{scope}

\begin{pgfonlayer}{background}

\filldraw[fill=black!50!, fill opacity=.2,rounded corners=20pt, line width=1pt, draw=black] (-7.5,-.5)--(-4,-.5)--(1.2,5)--cycle;

\filldraw[fill=black!50!, fill opacity=.2,rounded corners=20pt, line width=1pt, draw=black] (7.5,.5)--(4,.5)--(-1.2,-5)--cycle;

\filldraw[fill=black!50!, fill opacity=.2, rounded corners=5pt, line width=1pt, draw=black] (-1.8,1.8) rectangle (1.8,-1.8);
\end{pgfonlayer}

\path[-]
		(v1)
		edge[help lines] (w1)
		edge[help lines] (w2)
		edge[help lines] (w3)
		edge[help lines] (w4)
		(v2)
		edge[help lines] (u1)
		edge[help lines] (u2)
		edge[help lines] (u3)
		(v3)
		edge[help lines] (u1)
		edge[help lines] (u3)
		edge[help lines] (u4)
		(u1)
		edge[help lines] (u4)
		edge[help lines] (w1)
		edge[help lines] (w2)
		(u2)
		edge[help lines] (u3)
		edge[help lines] (w2)
		edge[help lines] (w3)
		(u3)
		edge[help lines] (w3)
		edge[help lines] (w4)
		(u4)
		edge[help lines] (w4)
		edge[help lines] (w1);

\node at (u3) [yshift=-1.5cm, font=\Large] {$\Delta(H)=6$}	;	

\end{tikzpicture}
%
}
\caption{A partition $(H_1,H_2)$ of $H$ such that $\text{col}(H_i) \leq 3$ and $\Delta(H_i) \leq 3$ for $i=1,2$.}
\end{figure}

The condition $\Delta(G) \geq 3$ in the simple case ensures that $G$ is not an odd cycle. However, since the hypergraphs of type (C) may have an arbitrary large maximum degree, we have to exclude this case manually. Before we prove Theorem~\sref{Theorem:col_partition}, it is necessary to obtain the following statement.

\begin{proposition}
If a hypergraph $H$ is $f$-partitionable for some $f \in \mathcal{V}_p(H)$ with $$f_1(v) + f_2(v) + \ldots + f_p(v) \geq d_H(v)$$ for all $v \in V(H)$, then there is an $f$-partition $(H_1,H_2,\ldots,H_p)$ of $H$ such that $d_{H_i}(v) \leq f_i(v)$ for all $v \in V(H_i)$ and for all $i \in \{1,2,\ldots,p\}$.
\end{proposition}

\begin{proof}
Given an arbitrary $p$-partition $(H_1,H_2,\ldots,H_p)$ of $H$, define its weight by
$$W_{(H_1,H_2,\ldots,H_p)}= \sum_{i=1}^p \big( |E(H_i)| - \sum_{v \in V(H_i)} f_i(v)  \big).$$
If there is a $v \in V(H)$ and two indices $i \neq j$  from $\{1,2,\ldots,p\}$  such that $v \in V(H_i)$, $d_{H_i}(v) \geq f_i(v)$ and $d_{H_j}(v) < f_j(v)$, then shifting $v$ from $H_i$ to $H_j$ decreases $W_{(H_1,H_2,\ldots,H_p)}$. In order to prove this, let $i < j$ (by symmetry), let $H_i'=H_i - v$, $H_j' = H_j + v$, and let $H_k'=H_k$ for all $k \in \{1,2,\ldots,p\} \setminus \{i,j\}$. Then, for $W=W_{(H_1,H_2,\ldots,H_p)}$ and for $W'=W_{(H_1',H_2',\ldots,H_p')}$ it holds
$$ W' - W = -d_{H_i} (v) + f_i(v) + d_{H_j} (v) - f_j(v) < 0.$$
Let $(H_1,H_2,\ldots,H_p)$ be an $f$-partition of $H$ that minimizes $W_{(H_1,H_2,\ldots,H_p)}$. We claim that $(H_1,H_2,\ldots,H_p)$ has the desired property. Otherwise there is be an index $i \in \{1,2,\ldots,p\}$ and a vertex $v \in V(H_i)$ such that $v \in V(H_i)$ and $d_{H_i}(v) > f_i(v)$. As $f_1(v) + f_2(v) + \ldots + f_p(v) \geq d_H(v)$, there is an index $j \in \{1,2,\ldots,p\}$ such that $d_{H_j}(v) < f_j(v)$. Thus, $H_j'=H_j + v$ is still strictly $f_j$-degenerate. Furthermore, $H_i'=H_i - v$ is strictly $f_i$-degenerate as well, and by the above observation, we obtain a new $p$-partition $(H_1',H_2',\ldots,H_p')$ with $W_{(H_1',H_2',\ldots,H_p')} < W_{(H_1,H_2,\ldots,H_p)}$, a contradiction.
\end{proof}

It is notable that the above proposition leads to a stronger version of Theorem~\ref{Theorem:Hauptsatz}.

\begin{stheorem}[\ref{Theorem:Hauptsatz}]
Let $H$ be a connected hypergraph, and let $f \in \mathcal{V}_p(H)$ be a vector function with $p \geq 1$ such that $f_1(v) + f_2(v) + \ldots + f_p(v) \geq d_H(v)$ for all $v \in V(H)$. Then, there is an $f$-partition $(H_1,H_2,\ldots,H_p)$ of $H$ such that $d_{H_i}(v) \leq f_i(v)$ for all $v \in V(H_i)$ and for all $i \in \{1,2,\ldots,p\}$ if and only if $(H,f)$ is not a hard pair.
\end{stheorem}

Now we are able to prove Theorem~\sref{Theorem:col_partition}.

\begin{proofof}[Theorem~\sref{Theorem:col_partition}]
Let $f_i(v) = k_i$ for all $v \in V(H)$ and for each $i \in \{1,2,\ldots,p\}$. Then, $f_1(v) + f_2(v) + \ldots + f_p(v) \geq \Delta(H) \geq d_H(v)$ and $f_i(v) \geq 1$ for all $i \in \{1,2,\ldots,p\}$ and for all $v \in V(H_i)$. Since $p \geq 2$, this implies that $(H,f)$ cannot be of type (M). Moreover, since $H$ is not a $tK_n$ for some $t,n \geq 1$ nor a $tC_n$ for $t\geq 1$ and $n \geq 5$ odd, it is easy to see that $(H,f)$ is not a hard pair (see Proposition~\ref{prop_block-hardpair}). Thus, by Theorem~\sref{Theorem:Hauptsatz}, $H$ admits an $f$-partition $(H_1,H_2,\ldots,H_p)$ such that $d_{H_i}(v) \leq f_i(v)=k_i$ for all $v \in V(H_i)$ and each $i \in \{1,2,\ldots,p\}$. In particular, $H_i$ is strictly $k_i$-degenerate and, thus, $\text{col}(H_i) \leq k_i$ for all $i \in \{1,2,\ldots,p\}$.
\end{proofof}

\subsection{(List)-point-partition number}
The \textbf{point-partition number} $\alpha_s(H)$ of a hypergraph $H$ (with $s \geq 0$) is the minimum number $k$ such that $H$ admits a $k$-coloring in which each color class induces an $s$-degenerate subhypergraph. Thus, $\alpha_0(H)$ corresponds to the chromatic number of $H$. Furthermore, the \textbf{list-point partition number} $\alpha_s^\ell(H)$ of a hypergraph $H$ is the least integer $k$ such that for any list-assignment $L$ fulfilling $|L(v)| \geq k$ for all $v \in V(H)$, there is an $L$-coloring of $H$ such that each color class induces an $s$-degenerate subhypergraph. For simple graphs, the point-partition number was introduced by Lick and White \cite{LiWi}. It is notable that for an arbitrary graph $G$, the \textbf{point arboricity} of $G$ is defined as the least number $k$ of forests forming a $k$-partition of $G$ and, thus, corresponds to $\alpha_1(G)$.

If we regard Theorem~\sref{Theorem:col_partition}, by setting $k_1=k_2=\ldots=k_p=s+1$, we obtain that the point-partition number $\alpha_s(H)$ is at most $p$ if $H$ is a connected hypergraph with maximum degree $\Delta \geq 1$ different from $tK_n$ with $t,n \geq 1$ and $tC_n$ for $t\geq1$ and $n \geq 3$ odd such that $p(s+1) =k_1 + k_2 + \ldots + k_p \geq \Delta$. For simple ordinary graphs, these cases were originally solved by Kronk and Mitchem \cite{KroMit} and Mitchem \cite{Mitch}. Concerning the list-point partition number, we obtain the following result as a simple consequence of Theorem~\ref{Theorem:Hauptsatz}.

\begin{corollary}\label{cor_list-point}
Let $H$ be a connected hypergraph different from $tK_n$ for $t,n \geq 1$ and different from $tC_n$ for $t \geq 1$, $n \geq 3$ odd, and let $\Delta(H) = \Delta \geq 1$. Furthermore, let $k$ and $s$ be integers such that $k \cdot s \geq \Delta$, $k \geq 2$, and let $L$ be a list-assignment such that $|L(v)| \geq k$ for all $v \in V(H)$. Then, there is an $L$-coloring of $H$ such that each color class induces a strictly $s$-degenerate subhypergraph.
\end{corollary}

\begin{proof}
Let $C= \bigcup_{v\in V(H)} L(v)$. By renaming the colors, we get $C=\{1,2,\ldots,p\}$. Let $f \in \mathcal{V}_p(H)$ be the function with
$$f_i(v)= \begin{cases}
s & \text{if } i \in L(v),\\
0 & \text{otherwise}
\end{cases}$$
for all $i \in \{1,2,\ldots,p\}$. Then,
\begin{align}
\label{sum_condition}
\sum_{i=1}^p f_i(v) \geq |L(v)|s \geq ks \geq \Delta \geq d_H(v)
\end{align}

for all $v \in V(H)$. If $H$ admits an $f$-partition $(H_1,H_2,\ldots,H_p)$, then setting $\varphi(v)=i$ if $v \in V(H_i)$ gives us the required $L$-coloring of $H$. Assume that $H$ is not $f$-partitionable. By \eqref{sum_condition} and Theorem~\ref{Theorem:Hauptsatz} it then follows that $(H,f)$ is a hard pair which implies, in particular, that $H$ is $\Delta$-regular and that $\Delta=ks$. Since $k \geq 2$, for each vertex $v$ there are two indices $i \neq j$ such that $f_i(v) \neq 0$ and $f_j(v) \neq 0$. Since $H$ is $\Delta$-regular, Proposition~\ref{prop_block-hardpair} implies that $H$ consists of just one block and $(H,f)$ is of type (K) or (C), a contradiction to the premise.
\end{proof}

Let $H$ be a hypergraph and let $L$ be an arbitrary list-assignment for $H$. We say that $H$ is $L \times s$\textbf{-choosable} if there is an $L$-coloring of $H$ such that each color class induces a strictly $s$-degenerate subhypergraph. With the help of this definition we can obtain a natural extension of Corollary~\ref{cor_list-point}. For simple ordinary graphs, it was proven by Borodin, Kostochka and Toft \cite{BorKosToft}. For hypergraphs, their proof can be copied as it stands.

\begin{theorem}
Let $H$ be a connected hypergraph. Then, $H$ is $L \times s$-choosable for each list-assignment $L$ satisfying $|L(v)| \geq d_H(v) /s$ for each $v \in V(H)$ if and only if at least one block of $H$ is different from $tK_n$ for all $t, n \geq 1$, from an $s$-regular hypergraph, and from a $tC_n$ with $t \geq 1$ and $n \geq 3$ odd.
\end{theorem}

\begin{proof}
As in the above corollary let $C$ be the set of colors used in the union of all lists $L(v)$ and assume $C=\{1,2,\ldots,p\}$. Moreover, define $f \in \mathcal{V}_p(H)$ with
$$f_i(v)=\begin{cases}
s & \text{if } i \in L(v),\\
0 & \text{otherwise}
\end{cases}$$ for all $i \in \{1,2,\ldots,p\}$. Then, $f_1(v) + f_2(v) + \ldots + f_p(v) \geq |L(v)|s \geq d_H(v)$ for all $v \in V(H)$ Clearly, $H$ is $L \times s$-choosable if and only if $H$ admits an $f$-partition. Note that $(H,f)$ is of type (M) if and only if $H$ is $s$-regular (by construction of the vector function $f$). Applying Theorem~\ref{Theorem:Hauptsatz} then completes the proof.
\end{proof}

\section{Concluding remarks}
It seems very likely that one can deduce a polynomial time algorithm from the proof of Theorem~\ref{Theorem:Hauptsatz}, which, given a hypergraph $H$ and a vector-function $f$ satisfying $f_1(v) + f_2(v) + \ldots + f_p(v) \geq d_H(v)$ for all $v \in V(H)$, finds an $f$-partition of $H$ or shows that $(H,f)$ is a hard pair. If we give up on the condition $f_1(v) + f_2(v) + \ldots + f_p(v) \geq d_H(v)$ for all $v \in V(H)$, then the decision problem whether $H$ admits an $f$-partition is {\NPC} since it contains the {\NPC} coloring problem. The complexity of the $f$-partition problem of simple graphs for constant vector-functions $f$ and $p=2$ was established in a series of papers \cite{AbuFH, WuYZ, YangY}. Let $r, s\geq 1$ be integers. A hypergraph $H$ is said to be \DF{$(r,s)$-partionable}, if there is a partition $(H_1,H_2)$ of $H$ such that $H_1$ is strictly $r$-degenerate and $H_2$ is strictly $s$-degenerate. The decision problem whether a simple graph with maximum degree $\Delta\geq 3$ admits an $(r,s)$-partition is polynomial time solvable if $\Delta=3$, or $r+s\geq \Delta$ or $(r,s)=(1,1)$; in all other cases the decicision problem is {\NPC}. Note that a graph $G$ is (1,1)-partitonable if and only if $G$ is bipartite.

For another application of Theorem~\ref{Theorem:Hauptsatz} regarding generalized hypergraph colorings we recommend taking a look at the forthcoming paper \cite{Schweser18}.

\end{document}